\RequirePackage{etoolbox}
\csdef{input@path}{{style/}{graphics/}}
\documentclass[aos,MSNbibl,seceqn,dvips]{arximspdf}
\usepackage{algorithm}
\usepackage{algpseudocode}
\usepackage{subeqn}
\usepackage{graphicx}

\doi{10.1214/14-AOS1283}
\volume{43}
\issue{2}
\pubyear{2015}
\firstpage{819}
\lastpage{846}
\docsubty{FLA}

\makeatletter
\newcommand{\lleft}{\left}
\newcommand{\rrvert}{\vert}
\newcommand{\rright}{\right}
\newcommand{\rrVert}{\Vert}
\newcommand{\llvert}{\vert}
\newcommand{\llVert}{\Vert}
\newcommand{\vvvert}{|\!|\!|}
\newcommand{\vvvertdu}{\big|\!\big|\!\big|}
\newcommand{\iint}{\int\!\!\!\int}
\newtheorem{theorem}{Theorem}
\newtheorem{lemma}{Lemma}
\newtheorem{proposition}{Proposition}
\newproclaim{definition}{Definition}
\newproclaim{example}{Example}
\makeatother

\begin{document}
\begin{frontmatter}

\title{The geometry of kernelized spectral clustering}
\runtitle{The geometry of spectral clustering}

\begin{aug}
\author[A]{\fnms{Geoffrey}~\snm{Schiebinger}\corref{}\thanksref{T1}\ead[label=e1]{gschiebinger@gmail.com}},
\author[A]{\fnms{Martin J.}~\snm{Wainwright}\thanksref{T2}\ead[label=e2]{wainwrig@stat.berkeley.edu}}
\and
\author[A]{\fnms{Bin}~\snm{Yu}\thanksref{T3}\ead[label=e3]{binyu@stat.berkeley.edu}}
\runauthor{G. Schiebinger, M. J. Wainwright and B. Yu}
\affiliation{University of California, Berkeley}
\address[A]{Department of Statistics\\
University of California, Berkeley\\
Berkeley, California 94720\\
USA\\
\printead{e1}\\
\phantom{E-mail: }\printead*{e2}\\
\phantom{E-mail: }\printead*{e3}}
\end{aug}
\thankstext{T1}{Supported by a Graduate Research Fellowship from the
NSF---Grant number DGE-1106400.}
\thankstext{T2}{Supported by NSF Grants CIF-31712-23800, and DMS-11-07000.}
\thankstext{T3}{Supported by NSF Grants DMS-11-07000, DMS-11-60319
(FRG), CDS\&E-MSS 1228246, ARO Grant W911NF-11-1-0114
and the Center for Science of Information (CSoI), an US NSF Science and
Technology Center, under Grant agreement CCF-0939370.}

%
\received{\smonth{4} \syear{2014}}
%
\revised{\smonth{10} \syear{2014}}

%
\begin{abstract}
Clustering of data sets is a standard problem in many areas of science
and engineering. The method of spectral clustering is based on
embedding the data set using a kernel function, and using the top
eigenvectors of the normalized Laplacian to recover the connected
components. We study the \mbox{performance} of spectral clustering in
recovering the latent labels of i.i.d. samples from a finite mixture of
nonparametric distributions. The difficulty of this label recovery
problem depends on the overlap between mixture components and how
easily a mixture component is divided into two nonoverlapping
components. When the overlap is small compared to the indivisibility of
the mixture components, the principal eigenspace of the
population-level normalized Laplacian operator is approximately spanned
by the square-root kernelized component densities. In the finite sample
setting, and under the same assumption, embedded samples from different
components are approximately orthogonal with high probability when the
sample size is large. As a corollary we control the fraction of samples
mislabeled by spectral clustering under finite mixtures with
nonparametric components.
\end{abstract}

%
\end{frontmatter}

\setcounter{footnote}{3}
\section{Introduction}

In the past decade, spectral methods have emerged as a powerful
collection of nonparametric tools for unsupervised learning, or
clustering. How can we recover information about the geometry or
topology of a distribution from its samples? Clustering algorithms
attempt to answer the most basic form of this question. One way in
which to understand
spectral clustering is as a relaxation of the NP-hard problem of
searching for the best graph-cut. Spectral
graph partitioning---using the eigenvectors of a matrix to find graph
cuts---originated in the early 1970s with the work of
Fiedler~\cite{fiedler} and of Donath and Hoffman~\cite{donathhoffman}.
Spectral clustering was introduced in machine learning, with
applications to clustering data sets and computing image segmentations
(e.g.,~\cite{shimalik,ng,shimeila}). The past decade has witnessed an
explosion of different spectral clustering algorithms. One point of
variation is that some use the eigenvectors of the kernel
matrix~\cite{daspec,kerneleca,donathhoffman}, or adjacency matrix in
the graph setting, whereas others use the eigenvectors of the
normalized Laplacian matrix~\cite{ng,shimeila,shimalik,fiedler}. This
division goes all the way back to the work of Donath and Hoffman, who
proposed using the adjacency matrix, and of Fiedler, who proposed
using the normalized Laplacian matrix.

In its modern and most popular form, the spectral clustering
algorithm~\cite{ng,shimalik} involves two steps: first, the
eigenvectors of the normalized Laplacian are used to embed the
dataset, and second, the $K$-means clustering algorithm is applied to
the embedded dataset. The normalized Laplacian embedding is an
attractive preprocessing step because the transformed clusters tend to
be linearly separable. Ng et~al. \cite{ng} show that, under certain
conditions on the empirical kernel matrix, an embedded dataset will
cluster tightly around well-separated points on the unit sphere. Their
results apply to a fixed dataset, and do not model the underlying
distribution of the data. Recently Yan et~al. \cite{yan} derived an
expression for the fraction of data misclustered by spectral clustering by
computing an analytical expression for the second eigenvector of the
Laplacian. They assumed that the similarity matrix is a small
perturbation away from the ideal block diagonal case.

The embedding defined by the normalized Laplacian has also been
studied in the context of manifold learning, where the primary focus
has been convergence of the underlying eigenvectors. This work is
motivated in part by the fact that spectral properties of the limiting
Laplace--Beltrami operator have long been known to shed light on the
connectivity of a manifold~\cite{cheeger}. The Laplacian eigenmaps of
Belkin and Niyogi~\cite{eigenmap} reconstruct Laplace--Beltrami eigenfunctions
from sampled data. Koltchinskii and
Gin\'{e}~\cite{kolgine} analyze the convergence of the empirical graph
Laplacian to the Laplace--Beltrami operator at a fixed point in the
manifold. von Luxburg and Belkin~\cite{luxburg} establish consistency
for the embedding in as much as the eigenvectors of the Laplacian
matrix converge uniformly to the eigenfunctions of the Laplacian
operator. Rosasco et~al. \cite{rosasco} provide simpler proofs of
this convergence, and in part, our work sharpens these results by
removing an unnecessary smoothness assumption on the kernel function.

In this paper, we study spectral clustering in the context of a
nonparametric mixture model. The study of spectral clustering under
nonparametric mixtures was initiated by Shi et~al. \cite{daspec}.
One of their theorems characterizes the top eigenfunction of a
kernel integral operator, showing that it does not change sign. One
difficulty in using the eigenfunctions of a kernel integral operator
to separate mixture components is that several of the top
eigenfunctions may correspond to a single mixture component (e.g., one
with a larger mixture weight). They propose that eigenfunctions of
the kernel integral operator that approximately do not change sign
correspond to different mixture components. However, their analysis
does not deal with finite datasets nor does it provide bounds on the
fraction of points misclustered.

The main contribution of this paper is an analysis of the normalized
Laplacian embedding of i.i.d. samples from a finite mixture with
nonparametric components. We begin by providing a novel and useful
characterization of the \mbox{principal} eigenspace of the population-level
normalized Laplacian operator: more precisely, when the mixture
components are indivisible and have small overlap, the eigenspace is
close to the span of the square root kernelized component densities.
We then use this characterization to analyze the geometric structure
of the embedding of a finite set of i.i.d. samples. Our main result
is to establish a certain geometric property of nonparametric mixtures
referred to as \emph{orthogonal cone structure}. In particular, we
show that
when the mixture components are indivisible and have small overlap,
embedded samples from different components are almost orthogonal with
high probability. We then prove that this geometric
structure allows $K$-means to correctly
label most of the samples. Our proofs rely on techniques from
operator perturbation theory, empirical process theory and spectral
graph theory.

The remainder of this paper is organized as follows. In
Section~\ref{SecBackground}, we set up the problem of separating the
components of a mixture distribution. We state our main results and
explore some of their consequences in Section~\ref{SecMain}. We prove
our main
results in Section~\ref{SecProofs}, deferring the proofs of several
supporting lemmas to the supplementary material \cite{SchiebingerWainwrightYu14Sup}.

\subsection*{Notation}
For a generic distribution
$\mathbb{P}$ on a measurable space ${\mathcal X}$, we denote the Hilbert
space of real-valued square integrable functions on ${\mathcal X}$ by
${L^2(\mathbb{P})}$. The ${L^2(\mathbb{P})}$ inner product is given by
${ \langle f,g \rangle}_{\mathbb{P}} = \int
f(x)g(x)\,d\mathbb{P}(x)$, and it induces the
norm $\llVert f \rrVert_{\mathbb{P}}$. The norm ${\llVert f \rrVert
_\infty}$ is the supremum of
the function $f$, up to sets of measure zero, where the relevant
measure is understood from context. The Hilbert--Schmidt norm of an
operator $\mathbf{T}\dvtx  {L^2(\mathbb{P})} \to{L^2(\mathbb{P})}$ is
$\vvvert \mathbf{T} \vvvert_{\mathrm{HS}}$, and the operator norm is $\vvvert \mathbf{T} \vvvert
_{\mathrm{op}}$.
The complement of a set
$B$ is denoted by $B^c$. See Appendix~D (supplementary material \cite{SchiebingerWainwrightYu14Sup}) for an
additional list of symbols.


\section{Background and problem set-up}\label{SecBackground}

We begin by introducing the family of nonparametric mixture models
analyzed in this paper, and then provide some background on kernel
functions, spectral clustering and Laplacian operators.

\subsection{Nonparametric mixture distributions}\label{SecMixture}

For some integer $K\geq2$, let
$\{\mathbb{P}_m\}_{m=1}^K$ be a collection of probability
measures on a compact space ${\mathcal X}$, and let the weights
$\{w_m\}_{m=1}^K$ belong to the relative
interior of the probability simplex in ${\mathbb{R}}^K$---that is,
$w_m\in(0,1)$ for all \mbox{$m= 1, \ldots,
K$,} and \mbox{$\sum_{m=1}^Kw_m=
1$.} This pair specifies a \emph{finite nonparametric mixture
distribution} via the convex combination
%
%
\begin{eqnarray}
\label{EqnDefnMixture} \bar{\mathbb{P}}&: =&\sum_{m= 1}^Kw_m
\mathbb{P}_m.
\end{eqnarray}
We refer to ${\{\mathbb{P}_m
\}_{m=1}^K}$ and ${\{w_m
\}_{m=1}^K}$ as the mixture components
and mixture weights, respectively. The family of
models~(\ref{EqnDefnMixture}) is \textit{nonparametric}, because the
mixture components are not constrained to any particular parametric
family.

A random variable ${\bar X}\sim\bar{\mathbb{P}}$ can be obtained
by first
drawing a multinomial random variable \mbox{${Z}\sim
\operatorname{Multinomial}(w_1, \ldots,
w_K)$,} and conditioning on the event \mbox{$\{ {Z}=
m\}$,} drawing a variable from mixture component
$\mathbb{P}_m$. Consequently, given a collection of samples
$\{X_i\}_{i=1}^{n}$ drawn i.i.d. from $\bar{\mathbb{P}}$, there is an
underlying set of \emph{latent labels} $\{ {Z}_i\}_{i=1}^{n}$.
Thus in the context of a mixture distribution, the clustering problem
can be formalized as recovering these latent labels based on observing
only the unlabeled samples $\{X_i\}_{i=1}^{n}$.

Of course, this clustering problem is ill defined whenever $\mathbb
{P}_j =
\mathbb{P}_k$ for some $j \neq k$. More generally, recovery of labels
becomes more difficult as the overlap of any pair $\mathbb{P}_j$ and
$\mathbb{P}_k$ increases, or if it is ``easy'' to divide any component into
two nonoverlapping distributions. This intuition is formalized in
our definition of the overlap and indivisibility parameters in
Section~\ref{SecKeyParameters} to follow.


\subsection{Kernels and spectral clustering}\label{SecKernels}

We now provide some background on spectral clustering methods and the
normalized Laplacian embedding. A kernel~$k$ associated with the
space ${\mathcal X}$ is a symmetric, continuous function
$k\dvtx  {\mathcal X}\times{\mathcal X}\to(0,\infty)$. A kernel is
said to be positive semidefinite if for any integer ${n}\geq1$
and elements $x_1,\ldots,x_{n}\in{\mathcal X}$, the kernel
matrix $A\in\mathbb{R}^{{n}\times{n}}$ with entries
\mbox{$A_{ij} = k(x_i,x_j)/n$} is
positive semidefinite. Throughout we consider a fixed but arbitrary
positive semidefinite kernel function.
In application to spectral clustering, one
purpose of a kernel function is to provide a measure of the
similarity between data points. A~canonical example is the Gaussian
kernel $k(x, x') = \exp(- \|x - x'\|_2^2)$; it
is close to $1$ for vectors $x$ and $x'$ that are relatively close, and
decays to zero for pairs that are far apart.

Let us now describe the normalized Laplacian embedding, which is a
standard part of many spectral clustering routines. Given ${n}$
i.i.d. samples $\{ X_i\}_{i=1}^{n}$ from $\bar{\mathbb{P}}$, the
associated kernel matrix $A\in{\mathbb{R}}^{{n}\times{n}}$
has entries $A_{ij} = \frac{1}{{n}}
k(X_i, X_j)$. The normalized Laplacian
matrix\footnote{To be precise, the matrix $I-L$ is actually the
normalized graph Laplacian matrix. However, the eigenvectors of
$L$ are identical to those of $I -L$, and we find it
simpler to work with $L$.} is obtained by rescaling the kernel
matrix by its row sums, namely
%
%
\begin{equation}
\label{EqnNormLapMat} L= D^{-1/2} AD^{-1/2},
\end{equation}
where\vspace*{1pt} $D$ is a diagonal matrix with entries $D_{ii} =
\sum_{j = 1}^{n}A_{ij}$. Since $L$ is a symmetric
matrix by construction, it has an orthonormal basis of eigenvectors,
and we let $\{v_1,\ldots,v_K\}$ denote the eigenvectors
corresponding to the largest $K$ eigenvalues of $L$. The
\emph{normalized Laplacian embedding} is defined on the basis of these
eigenvectors: it is the map
\mbox{$\Phi_{\mathcal V}\dvtx \{X_1,\ldots,X_{n}\} \to\mathbb{R}
^K$}
defined by
%
%
\begin{eqnarray}
\label{EqnSpectralEmbedding} \Phi_{\mathcal V}(X_i) &: =&(v_{1i},\ldots, v_{Ki}).
\end{eqnarray}
A typical form of spectral clustering consists of the following two
steps. First, compute the normalized Laplacian, and map each
data point $X_i$ to a $K$-vector via the
embedding~(\ref{EqnSpectralEmbedding}). The second step is to apply a
standard clustering method (such as $K$-means clustering) to the
embedded data points. The conventional rationale for the second step
is that the embedding step typically helps reveal cluster structure
in the data set, so that it can be found by a relatively simple
algorithm. The goal of this paper is to formalize the sense in which
the normalized Laplacian embedding~(\ref{EqnSpectralEmbedding}) has
this desirable property.


We do so by first analyzing the population operator that underlies the
normalized Laplacian matrix.
It is defined by the \emph{normalized kernel
function}
%
%
\begin{eqnarray}
\label{EqnDefnRescaledKernel} {\bar k}(x,y) &: =&\frac{1}{\bar{q}(x)}
k(x,y) \frac{1}{\bar{q}(y)},
\end{eqnarray}
where $\bar{q}(y) = \sqrt{\int k(x,y) \,d\bar{\mathbb{P}}(x)}$.
Note that this kernel
function can be seen as a continuous analog of the normalized
Laplacian matrix~(\ref{EqnNormLapMat}).

The normalized kernel function in conjunction with the mixture defines
the \emph{normalized Laplacian operator} $\bar{\mathbf{T}}\dvtx
L^2(\bar{\mathbb{P}})
\rightarrow L^2(\bar{\mathbb{P}})$ given by
%
%
\begin{eqnarray}
\label{EqnDefnNormLapOp} (\bar{\mathbf{T}}f) (\cdot) &: =&\int{\bar
k}(\cdot, y)f(y) \,d
\bar{\mathbb{P}}(y).
\end{eqnarray}
Under suitable regularity conditions (see Appendix~\textup{C.1} (supplementary\break material~\cite{SchiebingerWainwrightYu14Sup}) for
details), this operator has an orthonormal set of\break eigenfunctions---with
eigenvalues in $[0,1]$---and
our main results relate these eigenfunctions to the underlying mixture
components $\{{\mathbb{P}}_m\}_{m=1}^K$.


\section{Analysis of the normalized Laplacian embedding}\label{SecMain}

This section is devoted to the statement of our main results, and
discussion of their consequences. These results involve a few
parameters of the mixture distribution, including its overlap and
indivisibility parameters, which are defind in
Section~\ref{SecKeyParameters}. Our first main result
(Theorem~\ref{ThmPop} in Section~\ref{SecPop}) characterizes the
principal eigenspace of the population-level normalized Laplacian
operator~(\ref{EqnDefnNormLapOp}), showing that it approximately
spanned by the square root kernelized densities of the mixture
components, as defined in Section~\ref{SecKeyParameters}. Our second
main result (Theorem~\ref{ThmMain} in Section~\ref{SecFinite})
provides a quantitative description of the angular structure in the
normalized Laplacian embedding of a finite sample from a mixture
distribution.


\subsection{Cluster similarity, coupling and indivisibility parameters}\label{SecKeyParameters}

In this section, we define some parameters associated with any
nonparametric mixture distribution, as viewed through the lens of a
given kernel. These quantities play an important role in our main
results, as they reflect the intrinsic difficulty of the clustering
problem.

Our first parameter is the \emph{similarity index} of the mixture
components $\{{\mathbb{P}}_m\}_{m=1}^K$. For any pair of
distinct indices $\ell\neq m$, the ratio
\begin{eqnarray*}
{\mathcal{S}}({\mathbb{P}}_\ell, {\mathbb{P}}_m) &: =&
\frac
{\int_{{\mathcal X}}
\int_{{\mathcal X}} k(x,y) \,d {\mathbb{P}}_m(x) \,d {\mathbb{P}}_\ell
(y)} {
\int_{{\mathcal X}} \int_{{\mathcal X}} k(x,y) \,d {{\bar{\mathbb
{P}}}}(x) \,d
{\mathbb{P}}_\ell(y)}
\end{eqnarray*}
is a kernel-dependent measure of the expected similarity between the
clusters indexed
by ${\mathbb{P}}_\ell$ and ${\mathbb{P}}_m$, respectively. Note that
${\mathcal{S}}
$ is not symmetric in
its arguments. The \emph{maximum similarity} over all
mixture components
%
%
\begin{eqnarray}
\label{EqnDefnSIMMAX} \mathcal{S}_{\max}({\bar{\mathbb{P}}}) &: =&\mathop{\max
_{\ell, m= 1, \ldots,K}}_{\ell\neq m} {\mathcal{S}}({\mathbb{P}}_\ell,
{\mathbb{P}}_m)
\end{eqnarray}
measures the overlap between mixture components with respect to the
kernel $k$.

Our second parameter, known as the coupling~parameter, is defined in
terms of
the square root kernelized densities of the mixture components. More
precisely, given any distribution $\mathbb{P}$, its \emph{square root
kernelized density} is the function \mbox{$q\in{L^2(\mathbb{P})}$}
given by
%
%
\begin{eqnarray}
\label{EqnDefnSquareRoot} q(x) &: =&\sqrt{\int k(x,y) \,d\mathbb{P}(y)}.
\end{eqnarray}
In particular, we denote the square root kernelized density of the
mixture distribution $\bar{\mathbb{P}}$ by $\bar{q}$, and those of
the mixture
components $\{\mathbb{P}_m\}_{m=1}^K$ by
$\{q_m\}_{m=1}^K$. In analogy with the normalized kernel
function ${\bar k}$ from equation~(\ref{EqnDefnRescaledKernel}), we also
define a normalized kernel for each mixture component, namely
%
%
\begin{eqnarray}
\label{EqnMixtureCompRescaled} k_m(x,y) &: =&\frac
{k(x,y)}{q_m(x)q_m(y)} \qquad\mbox{for $m=
1, \ldots, K$.}
\end{eqnarray}
The \emph{coupling parameter}
%
%
\begin{eqnarray}
\label{EqnDefnKERGAP} {{\mathcal{C}}({\bar{\mathbb{P}}})} &: =&\max
_{m= 1, \ldots, K} {\llVert{k_m - w_m{\bar k}}
\rrVert^2_{\mathbb{P}_m\otimes\mathbb{P}_m}}
\end{eqnarray}
measures the coupling of the spaces ${L^2(\mathbb{P}_m)}$ with
respect to $\bar{\mathbf{T}}$. In particular, when ${{\mathcal
{C}}({\bar{\mathbb{P}}})} = 0$, then
the normalized Laplacian can be decomposed as the sum
%
%
\begin{eqnarray}
\label{EqnLapDecomp} \bar{\mathbf{T}}& =& \sum_{m= 1}^Kw_m
\mathbf{T} _m,
\end{eqnarray}
where $(\mathbf{T}_mf)(y) = \int f(x) k_m(x,y) \,d \mathbb{P}_m
(x)$ is
the operator defined by the normalized kernel $k_m$. When the
coupling parameter is no longer exactly zero but still small, then
decomposition~(\ref{EqnLapDecomp}) still holds in an approximate
sense.

%

Our final parameter measures how easy or difficult it is to ``split''
any given mixture component $\mathbb{P}_m$ into two or more parts. If
this splitting can be done easily for any component, then the mixture
distribution will be hard to identify, since there is an ambiguity as
to whether $\mathbb{P}_m$ defines one component or multiple
components. In order to formalize this intuition, for a distribution
$\mathbb{P}$ and for a measurable subset $S\subset{\mathcal X}$, we
introduce the shorthand notation \mbox{$p(S) = \int_S
\int_{\mathcal X}k(x,y) \,d\mathbb{P}(x)\,d\mathbb{P}(y)$}. With this
notation, the
indivisibility of $\mathbb{P}$ is
%
%
\begin{eqnarray}
\label{EqnDefnIndivisiblity} {\Gamma({\mathbb{P}})} &: =&\inf_S
\frac{ p({\mathcal X})
\int_S
\int_{S^c} k(x,y) \,d\mathbb{P}(x) \,d\mathbb{P}(y)}{p(S)p(S^c)},
\end{eqnarray}
where the infimum is taken over all measurable subsets $S$ such that
$p(S) \in(0,1)$. The \emph{indivisibility parameter}
$\Gamma_{\min}({\bar{\mathbb{P}}})$ of a mixture distribution
${\bar{\mathbb{P}}}$ is the minimum
indivisibility of its mixture components
%
%
\begin{eqnarray}
\label{EqnDefnCHEEGMIN} \Gamma_{\min}({\bar{\mathbb{P}}}) &: =&\min
_{m= 1, \ldots, K} {\Gamma({\mathbb{P}_m})}.
\end{eqnarray}

Our results in the next section apply when the similarity ${\mathcal
{S}}_{\max}({\bar{\mathbb{P}}})$ and coupling ${{\mathcal
{C}}({\bar{\mathbb{P}}})}$ are small compared to the
indivisibility $\Gamma_{\min}({\bar{\mathbb{P}}})$.
Some examples help illustrate when this is the case.

%
\begin{example}\label{ExTriangle}
Consider the one-dimensional triangular density function
\[
g_{{\mathbb{T}}_{\mu}}(x):= \cases{ x - \mu+ 1, &\quad if $x \in(\mu-
1, \mu)$;
\cr
-x
+ \mu+ 1, &\quad if $x \in(\mu, \mu+ 1)$;
\cr
0, &\quad otherwise,}
\]
with location $\mu> 0$. We denote corresponding distribution by
${\mathbb{T}}_{\mu}$. In this example we calculate the similarity,
coupling and indivisibility parameters for the mixture of triangular
distributions $\bar{\mathbb{T}}: =\frac{1}{2}{\mathbb{T}}_{0} +
\frac{1}{2}{\mathbb{T}}_{\mu}$ and the uniform kernel
$k_\nu(x,y) = \frac{1}{2\nu}\mathbf{1}\{\llvert x-y \rrvert\le
\nu\}$
with bandwidth $\nu\in(0,1)$.\footnote{This is
not a positive semidefinite kernel function, but it helps to build
intuition our intuition in a case where all the integrals are easy.}

\subsubsection*{Similarity}
It is
straightforward to calculate the similarity parameter ${\mathcal
{S}}_{\max}(\bar{\mathbb{T}})$
by solving a few simple integrals. We find that
\[
\mathcal{S}_{\max}(\bar{\mathbb{T}}) = \frac{2(2 + \nu- \mu
)_+^4}{\nu(16 -
8\nu^2 + 3\nu^3) + (2 + \nu- \mu)_+^4}.
\]

\subsubsection*{Coupling}
To compute the~coupling~parameter, we must compute the kernelized
densities of ${\mathbb{T}}_{0}$ and $ {\mathbb{T}}_{\mu}$, and the
normalized kernel functions $k_1(x,y)$ and ${\bar k}(x,y)$. Some
calculation yields the following equation for the kernelized density:
%
%
\begin{equation}
\label{EqnTriangKerDens} q_1^2(x) =
\cases{
\displaystyle\frac{(1 + \nu- x)^2}{4\nu}, &\quad if $x \in(-1 - \nu, -1 + \nu)$,
\vspace*{5pt}\cr
\displaystyle 1 - \frac{\nu}{2} - \frac{x^2}{2\nu}, &\quad if $x \in(-\nu,\nu)$,
\vspace*{5pt}\cr
\displaystyle\frac{(1 + \nu+ x)^2}{4\nu}, &\quad if $x \in(1 - \nu, 1 + \nu)$,
\vspace*{5pt}\cr
\displaystyle g_{{\mathbb{T}}_{0}}(x), & \quad otherwise.}
\end{equation}
As can be seen in Figure~\ref{FigTriangKerDens}, the kernelized
density $q_1^2(x)$ of ${\mathbb{T}}_{0}$ is a smoothed version of
$g_{{\mathbb{T}}_{0}}(x)$ that interpolates quadratically around the
nondifferentiable points of $g_{{\mathbb{T}}_{0}}(x)$.

%
\begin{figure}

\includegraphics{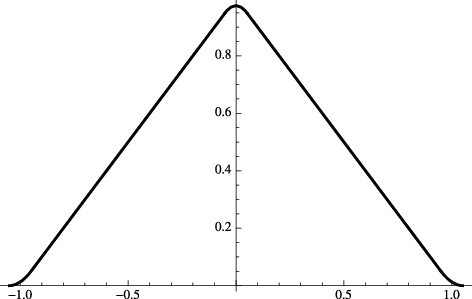}

\caption{The kernelized density $q_1^2(x)$ of equation~(\protect\ref
{EqnTriangKerDens}) with $\nu= 0.05$.}
\label{FigTriangKerDens}
\end{figure}

The kernelized density of ${\mathbb{T}}_{\mu}$ has the same shape
as that of ${\mathbb{T}}_{0}$ but is shifted by~$\mu$. In
particular,
\[
q_2^2(x) = q_1^2(x - \mu).
\]
Therefore the normalized kernels satisfy $k_1(x,y) = \frac{1}{2}
{\bar k}(x,y)$ for $x,y \in(-1, \mu- 1 - \nu)$. By upper
bounding the integrand over the remaining region, we find that
\[
{{\mathcal{C}}(\bar{\mathbb{T}})}^2 = \iint(k_1 - {{
\bar k}}/{2})^2 \,d\mathbb{P}_1(x,y) \le{2(2 + \nu-
\mu)_+}.
\]

\subsubsection*{Indivisibility}
It is straightforward to calculate the indivisibility
$\Gamma({\mathbb{T}}_{\mu})$. For any $\nu\in
(0,1)$ and $\mu\in\mathbb{R}$, the set $S$ defining
${\Gamma({{\mathbb{T}}_{\mu}})}$ is $S = (\mu,\infty)$. Hence
by solving a few simple integrals, we find that
\[
{\Gamma({{\mathbb{T}}_{\mu}})} = \frac{2(6-\nu
)(2-\nu)\nu}{16 - 8\nu^2 + 3\nu^3}.
\]
Note that ${\Gamma({{\mathbb{T}}_{\mu}})}$ does not depend on
$\mu$. Therefore the indivisibility of $\bar{\mathbb{T}}: =
\frac{1}{2} {\mathbb{T}}_{0} + \frac{1}{2}{\mathbb{T}}_{\mu}$ is
\[
\Gamma_{\min}(\bar{\mathbb{T}}) = {\Gamma({{\mathbb{T}}_{\mu}})}
= {\Gamma({{\mathbb{T}}_{0}})}.
\]
It is instructive to consider the indivisibility of the following
poorly defined two-component mixture:
\[
\bar{\mathbb{T}}_{\mathrm{bad}}: =\tfrac{1}{2}
\mathbb{P}_{\mathrm{bad}}(
\mu) + \tfrac{1}{2} {\mathbb{T}}_{2\mu},
\]
where $\mathbb{P}_{\mathrm{bad}}(\mu)$ is
the bimodal component $\mathbb{P}_{\mathrm{bad}}(\mu): =\frac{1}{2} {\mathbb{T}}_{0} +
\frac{1}{2}{\mathbb{T}}_{\mu}$. It is easy to verify that for any
$\nu\in(0,1)$ and $\mu> 2 + \nu$, the
indivisibility of the bimodal component is
$\Gamma(\mathbb{P}_{\mathrm{bad}}(\mu))
= 0$, and therefore
$\Gamma_{\min}(\bar{\mathbb{T}}_{\mathrm{bad}}) = 0$.
\end{example}

From Example~\ref{ExTriangle} we learn that the similarity ${\mathcal
{S}}_{\max}(\bar{\mathbb{T}})$ and
coupling ${{\mathcal{C}}(\bar{\mathbb{T}})}$ parameters decrease as
the offset $\mu$
increases. Together, these two parameters measure the overlap which
our intuition tells us should decrease as $\mu$ increases. On the
other hand, the indivisibility parameter $\Gamma_{\min}(\bar{\mathbb
{T}})$ is
independent of $\mu$.


Our next example is more realistic in the sense that the kernel and
mixture components do not have bounded support, and the kernel
function is positive semidefinite.
%

\begin{example}\label{ExGaussian}
In this example we calculate the similarity, coupling and
indivisibility parameters for the mixture of Gaussians $\bar\mathbb
{N}= \frac
{1}{2} \mathbb{N}(0,1) + \frac{1}{2} \mathbb{N}(\mu,1)$
equipped with the Gaussian kernel
$k_\nu(x,y) = \frac{1}{\sqrt{2\pi}\nu} \exp
[-\frac{\llvert x-y \rrvert^2}{2\nu^2} ]$.

\subsubsection*{Similarity}
As in the previous example, it is straightforward to calculate the
maximal intercluster similarity $\mathcal{S}_{\max}(\bar\mathbb
{N})$ by solving a handful
of Gaussian integrals. We find that
%
%
\begin{equation}
\mathcal{S}_{\max}(\bar\mathbb{N}) = \frac{2\exp({{-\mu
^2}/(2\nu^2 +4)})}{1 + \exp({{-\mu^2}/(2\nu^2 + 4)})}
\le4 e^{{-\mu^2}/(2\nu^2 + 4)}.
\end{equation}

\subsubsection*{Coupling}
The kernelized density of $\mathbb{N}(0,1)$ is
\[
q_1^2(x) = \frac{1}{\sqrt{\nu^2+1}}\exp\biggl[{
\frac
{-x^2}{2(\nu^2+1)}} \biggr],
\]
and the kernelized density of $\mathbb{N}(\mu,1)$ is simply the
translation $q_2^2(x) = q_1^2(x-\mu)$.
We can bound the coupling parameter ${{\mathcal{C}}(\bar\mathbb
{N})}$ by upper
bounding the integrand $k_1 - {\bar k}$ over a high-probability compact
set (a modification of the trick from Example~\ref{ExTriangle}).
We show the resulting bound in Figure~\ref{FigCouplingSim}. This
(albeit loose) bound captures the exponential decay of ${{\mathcal
{C}}(\bar\mathbb{N})}$ with $\mu$.

%
\begin{figure}

\includegraphics{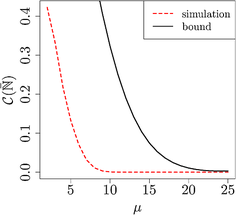}

\caption{The coupling parameter ${{\mathcal{C}}(\bar\mathbb{N})}$
for the mixture
of Gaussians with Gaussian kernel and $\nu= 2$. The red line
displays the simulated value of ${{\mathcal{C}}(\bar\mathbb{N})}$ as
a function of
the offset $\mu$ between mixture components. The black line
displays our analytical bound.}
\label{FigCouplingSim}
\end{figure}

\subsubsection*{Indivisibility}
It is straightforward to compute the indivisibility of the
unit-variance normal distribution $\mathbb{N}(\mu,1)$ with
location $\mu\in\mathbb{R}$. The set defining the indivisibility is $S
= (\mu,\infty)$. Solving a handful of Gaussian integrals yields
%
%
\begin{equation}
\label{EqnGaussianIndivisibility} \Gamma_{\min}\bigl(\mathbb{N}(\mu,1)\bigr) =
\frac{2}{\pi}\arctan\bigl({\nu\sqrt{2 + \nu^2}}\bigr).
\end{equation}
%
\end{example}

We conclude with a counter example, showing that the similarity
parameter is \emph{not} relatively small for the linear kernel $k(x,y)
= x\cdot y$.

%
%
\begin{example}\label{ex3}
Consider the mixture ${\bar{\mathbb{P}}}= \frac{1} 2 \mathbb{P}_1 +
\frac{1} 2 \mathbb{P}_2$
with components $\mathbb{P}_1$ uniform over $(1,2)$ and $\mathbb
{P}_2$ uniform
over $(2+\delta,3+\delta)$. With the linear kernel $k(x, y) = x \cdot
y$, we have
\[
\mathcal{S}_{\max}({\bar{\mathbb{P}}}) = \frac{\iint x y \,d\mathbb
{P}_1(x) \,d\mathbb{P}_2(y)}{\iint x y
\,d\mathbb{P}_1(x) \,d{\bar{\mathbb{P}}}(y)} =
\frac{\int y \,d\mathbb
{P}_2(y)}{\int y\,d{\bar{\mathbb{P}}}(y)} \ge\frac{1}{2}.
\]
Since $\Gamma_{\min}({\bar{\mathbb{P}}})$ is always between $0$
and $1$, this
calculation demonstrates that the similarity parameter
$\mathcal{S}_{\max}({\bar{\mathbb{P}}})$ is never small
compared to $\Gamma_{\min}({\bar{\mathbb{P}}})$.
\end{example}


\subsection{Population-level analysis}\label{SecPop}

In this section, we present our population-level analysis of the
normalized Laplacian embedding. Consider the following two subspaces
of ${L^2(\bar{\mathbb{P}})}$:
\begin{itemize}
\item the subspace $\mathcal{R}\subset{L^2(\bar{\mathbb{P}})}$
spanned by the
top $K$ eigenfunctions of the normalized Laplacian operator
$\bar{\mathbf{T}}$ from equation~(\ref{EqnDefnNormLapOp}) and
\item the span ${\mathcal Q}= \operatorname{span}\{q_1,\ldots,q_K\}
\subset{L^2(\bar{\mathbb{P}})}$ of the square root kernelized
densities; see
equation~(\ref{EqnDefnSquareRoot}).
\end{itemize}
The subspace ${\mathcal Q}$ can be used to define a map
$\Phi_{{\mathcal Q}}\dvtx  {\mathcal X}\to\mathbb{R}^K$ known as the
\emph{square-root kernelized density embedding}, given by
%
%
\begin{eqnarray}
\label{Eqkerdensembed} \Phi_{{\mathcal Q}}(x) &: =&\bigl(q_1(x), \ldots,q_K(x)\bigr).
\end{eqnarray}
This map is relevant to clustering, since the vector
$\Phi_{{\mathcal Q}}(x)$ encodes sufficient information to perform a
likelihood ratio test (based on the kernelized densities) for labeling
data points.

On the other hand, the subspace $\mathcal{R}$ is important because it
is the\break population-level quantity that underlies spectral clustering.
As described in Section~\ref{SecKernels}, the first step of spectral
clustering involves embedding the data using the eigenvectors of the
Laplacian matrix. This\vspace*{1pt} procedure can be understood as a way of
estimating the \emph{population-level Laplacian embedding}: more
precisely, the map $\Phi_{\mathcal{R}}\dvtx {\mathcal X}\to\mathbb{R}^K$ given
by
%
%
\begin{eqnarray}
\Phi_{\mathcal{R}}(x) &: =&\bigl(r_1(x), \ldots, r_K(x)
\bigr),
\end{eqnarray}
where $ \{r_m\}_{m=1}^K$ are the top $K$
eigenfunctions of the kernel operator $\bar{\mathbf{T}}$.

To build intuition, imagine for the moment varying the kernel function
so that the kernelized densities converge to the true densities. For
example, imagine sending the bandwidth of a Gaussian kernel to $0$.
While the kernelized densities \mbox{approach} the true densities, the
subspace $\mathcal{R}$ is only a well defined mathematical object for
kernels with nonzero bandwidth. Indeed, as the bandwidth shrinks to
zero, the eigengap separating the principal eigenspace $\mathcal{R}$ of
$\bar{\mathbf{T}}$ from its lower eigenspaces vanishes. For this
reason, we
analyze an arbitrary but fixed kernel function, and we discuss kernel
selection in Section~\ref{SecDisc}.

The goal of this section is to quantify the difference between the two
mappings $\Phi_{{\mathcal Q}}$ and $\Phi_{\mathcal{R}}$, or
equivalently between
the underlying subspaces ${\mathcal Q}$ and $\mathcal{R}$. We assume that
the square root kernelized densities $q_1,\ldots,q_K$ are linearly
independent so that ${\mathcal Q}$ has the same dimension, $K$, as
$\mathcal{R}$. This condition is very mild when the overlap parameters
$\mathcal{S}_{\max}(\bar{\mathbb{P}})$ and ${{\mathcal
{C}}(\bar{\mathbb{P}})}$ are small. We measure the
distance between these subspaces by the Hilbert--Schmidt
norm\footnote{Recall that the Hilbert--Schmidt norm of an operator is
the infinite dimensional analogue of the Frobenius norm of a
matrix.} applied to the difference between their orthogonal
projection operators,
%
%
\begin{eqnarray}
\label{EqnDefnDistance} \rho({\mathcal Q},\mathcal{R}) &: =&\vvvert
{\Pi_{{\mathcal
Q}}}
- {\Pi_{\mathcal{R}}} \vvvert_{\mathrm{HS}}.
\end{eqnarray}

Recall\vspace*{1pt} the similarity parameter $\mathcal{S}_{\max}({{\bar{\mathbb
{P}}}})$, coupling parameter
${{\mathcal{C}}({\bar{\mathbb{P}}})}$ and indivisibility parameter
$\Gamma_{\min}({\bar{\mathbb{P}}})$,
as previously defined in equations~(\ref{EqnDefnSIMMAX}), (\ref{EqnDefnKERGAP}) and~(\ref{EqnDefnCHEEGMIN}), respectively. Our
main results involve a function of these three parameters and the
minimum \mbox{$w_{\min}: =\min_{m= 1, \ldots, K}
w_m$} of the mixture weights, given by
%
%
\begin{eqnarray}
\label{EqnDifficulty} \varphi({\bar{\mathbb{P}}}; k) &: =&\frac{ \sqrt
{K}[\mathcal{S}_{\max}({\bar{\mathbb{P}}}) +
{{\mathcal{C}}({\bar{\mathbb{P}}})}]^{1/2}}{w_{\min} \Gamma
^2_{\min}({\bar{\mathbb{P}}})}.
\end{eqnarray}
Our first main theorem guarantees that as long as the mixture is
relatively well separated, as measured by the \emph{difficulty
function} $\varphi$, then the
$\rho$-distance~(\ref{EqnDefnDistance}) between $\mathcal{R}$ and
${\mathcal Q}$ is proportional to $\varphi({\bar{\mathbb{P}}};k)$.
Our theorem also
involves the quantity
\[
b_{\max}: =\max\limits
_{m= 1, \ldots, K} {\biggl\llVert\int k_m(x,y)\,d
\mathbb{P}_m(y) \biggr\rrVert_\infty^2}.
\]
Note that this is simply a constant whenever the kernels $k_m$ are
bounded.
%

\begin{theorem}[(Population control of subspaces)]
\label{ThmPop}
For any\vspace*{1pt} finite mixture ${\bar{\mathbb{P}}}$ with difficulty
function bounded as
\mbox{$\varphi({\bar{\mathbb{P}}};k) \le[576 \sqrt{12 +
b_{\max}}
]^{-1}\Gamma^2_{\min}({\bar{\mathbb{P}}}) $,} the distance
between subspaces
${\mathcal Q}$ and $\mathcal{R}$ is bounded as
%
%
\begin{eqnarray}
\label{EqnPopBound} \rho({\mathcal Q},\mathcal{R}) & \leq& 16 \sqrt{12 +
{b_{\max}}} \varphi({\bar{\mathbb{P}}};k).
\end{eqnarray}
\end{theorem}

%
\begin{figure}[t]

\includegraphics{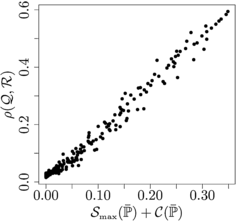}

\caption{The $\rho$-distance between ${\mathcal Q}$ and $\mathcal{R}$
scales linearly with $\mathcal{S}_{\max}({\bar{\mathbb{P}}}) +
{{\mathcal{C}}({\bar{\mathbb{P}}})}$. Each
point corresponds to a different offset $\mu$ between the two
Gaussian mixture components from Example~\protect\ref{ex3}.}\label{FigSubspaceError}
\end{figure}

Relationship~(\ref{EqnPopBound}) is easy to understand in the
context of translated copies of identical mixture components.
Consider the mixture of Gaussians with Gaussian kernel setup in
Example~\ref{ex3}. Recall from equation~(\ref{EqnGaussianIndivisibility})
that the indivisibility parameter is independent of the offset
$\mu$. Hence in this setting relationship~(\ref{EqnPopBound})
simplifies to
\[
\rho({\mathcal Q},\mathcal{R}) \le c\bigl[\mathcal{S}_{\max}({{\bar{
\mathbb{P}}}}) + {{\mathcal{C}}({\bar{\mathbb{P}}})}\bigr]^{1/2}.
\]
Figure~\ref{FigSubspaceError} shows a clear linear relationship
between $\rho({\mathcal Q},\mathcal{R})$ and ${\mathcal{S}}_{\max
}({\bar{\mathbb{P}}}) +
{{\mathcal{C}}({\bar{\mathbb{P}}})}$, suggesting that it might be
possible to remove the
square root in the clustering difficulty~(\ref{EqnDifficulty}).

\begin{figure}[b]

\includegraphics{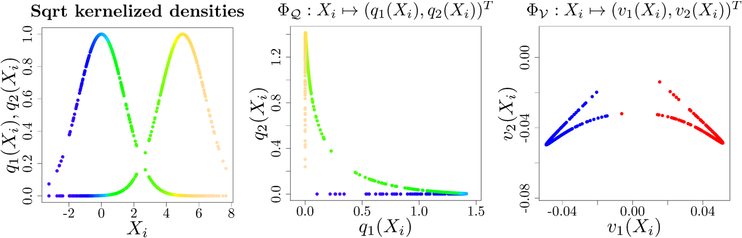}

\caption{Geometric structure in the square root kernelized density
embedding. (Left) square root kernelized densities for a mixture of
Gaussians with Gaussian kernel. The color of the $i$th dot
indicates the likelihood ratio of the mixture components at
$X_i$. (Center) data embedded under the square root kernelized
density embedding $\Phi_{{\mathcal Q}}$, colored by likelihood ratio.
(Right) normalized Laplacian embedding of the samples, colored by
latent label.} 
\label{FigKerDens}
\end{figure}

One important consequence of relationship~(\ref{EqnPopBound}) stems
from geometric structure in the
square root kernelized density embedding.
When there is little overlap
between mixture components with respect to the kernel, the square root
kernelized densities are
not simultaneously large; that is, $\Phi_{{\mathcal Q}}(X)$ will have at
most one component much different from zero. Therefore the data will
concentrate in tight spikes about the axes. This is illustrated in
Figure~\ref{FigKerDens}.


%
%


\subsection{Finite sample analysis}\label{SecFinite}

Thus far, our analysis has been limited to the population level,
corresponding to the ideal case of infinitely many samples. We now
turn to the case of finite samples. Here an additional level of
analysis is required, in order to relate empirical versions (based on
the finite collection of samples) to their population analogues.
Doing so allows us to show that under suitable conditions, the
Laplacian embedding applied to i.i.d. samples drawn from a finite
mixture satisfies a certain geometric property, which we call
\emph{orthogonal cone structure}, or OCS for short.

We begin by providing a precise definition of when an embedding
$\Phi\dvtx  {\mathcal{X}}\rightarrow{\mathbb{R}}^K$ reveals orthogonal cone
structure. Given a collection of labeled samples $\{X_i,
Z_i\}_{i=1}^{n}\subset{\mathcal{X}}\times[K] $ drawn from a
$K$-component
mixture distribution, we let ${\mathcal{Z}}_m= \{i \in[{n}]
|
Z_i = m\}$ denote the subset of samples drawn from mixture
component $m= 1, \ldots, K$. For any set ${\mathcal{Z}}\subseteq
[{n}] = \{1, 2, \ldots, {n}\}$, we\vspace*{1pt} use $|{\mathcal{Z}}|$ to denote its
cardinality. For vectors $u, v \in{\mathbb{R}}^{n}$, we use
$\operatorname{angle}(u, v) = \mbox{arccos} \frac{{\langle u, v
\rangle}}{\|u\|_2
\|v\|_2}$ to\vspace*{1pt} denote the angle between them. With this notation, we
have the following:

%
\begin{definition}[{[Orthogonal cone structure (OCS)]}]\label{DefAngStruct}
Given parameters $\alpha\in(0,1)$ and $\theta\in(0,\frac{\pi}4)$, the
embedded data set $\{\Phi(X_i), Z_i\}_{i=1}^{n}$ has
$(\alpha,\theta)$-OCS if there is an orthogonal basis
$\{e_1, \ldots, e_K\}$ of ${\mathbb{R}}^K$ such
that
\begin{eqnarray*}
\bigl\llvert\bigl\{ i \in[{n}] |\operatorname{angle} \bigl(\Phi(X_i),
e_m \bigr) < \theta\bigr\} \cap{\mathcal{Z}}_m \bigr
\rrvert& \geq& (1-\alpha) \llvert{\mathcal{Z}}_m\rrvert\qquad
\mbox{for all $m= 1,\ldots,K$.}
\end{eqnarray*}
\end{definition}

\begin{figure}[b]

\includegraphics{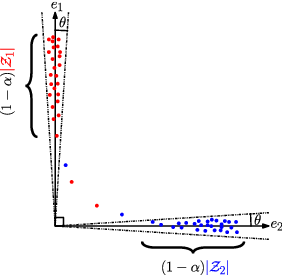}

\caption{Visualizing $(\alpha,\theta)$-OCS: the
labeled set of points plotted above has $(\alpha,\theta)$
orthogonal cone structure with respect to its labeling. The color of each
dot indicates the value of the corresponding label ${Z}_i\in
\{1,2\}$, where $1$ corresponds to red and $2$ to blue. This set of
points has $(\alpha,\theta)$-orthogonal cone structure because a
fraction $1-\alpha$ of the red points (for which ${Z}_i=
1$) lie with an angle $\theta$ of $e_1$, a fraction
$1-\alpha$ blue points (for which ${Z}_i= 2$) lie with
an angle $\theta$ of $e_2$ and $e_1$ is orthogonal to~$e_2$.}
\label{FigConeStruct}
\end{figure}

In words, a labeled dataset has orthogonal cone structure if
most pairs of embedded data points with distinct labels are almost
orthogonal. See Figure~\ref{FigConeStruct} for an illustration of
this property.

%

Our main theorem in the finite sample setting establishes that under
suitable conditions, the normalized Laplacian embedding has orthogonal
cone structure. In order to state this result precisely, we require a
few additional conditions.

\subsubsection*{Kernel parameters}
As a consequence of the compactness of ${\mathcal X}$, the kernel
function is $b$-bounded, meaning
that {$k(x, x') \in(0, b)$} for all $x, x' \in{\mathcal{X}}$. As
another consequence, the kernelized densities are lower bounded as
$q_m(X^m)\geq r>0$ with ${\bar{\mathbb{P}}}$-probability one.
In the following
statements, we use $c, {c}_0, {c}_1, \ldots$ to denote quantities
that may depend on $b$, and $r$ but are otherwise
independent of the mixture distribution.

\subsubsection*{Tail decay}
The tail decay of the mixture components enters our finite sample
result through the function $\psi\dvtx (0,\infty)\to[0,1]$,
defined by
%
%
\begin{eqnarray}
\label{EqnDefnTail} \psi(t) &: =&\sum_{m=1}^K{
\mathbb{P}}_m \biggl[\frac{q_m^2(X)}{\llVert q_m \rrVert_{{\bar
{\mathbb
{P}}}}^2} < t \biggr].
\end{eqnarray}
Note that $\psi$ is an increasing function with $\psi(0) = 0$.
The rate of increase of $\psi$ roughly
measures the tail decay of the square root kernelized densities.
Intuitively, perturbations to the square root kernelized density
embedding will have a greater effect on points closer to the origin.

Recall the population level clustering difficulty parameter
$\varphi({\bar{\mathbb{P}}};k)$ previously defined in
equation~(\ref{EqnDifficulty}).
Our theory requires that there is some $\delta> 0$ such that
%
%
\begin{eqnarray}
\label{EqnDevCondition} \underbrace{ \biggl[\varphi({\bar{\mathbb{P}}};
k) +
\frac{1}{\Gamma^2_{\min}({\bar{\mathbb{P}}})} \biggl(\frac
{1}{\sqrt{{n}}} + \delta\biggr)
\biggr]}_{\varphi_n(\delta)} & \leq& c \Gamma^2_{\min}({\bar{
\mathbb{P}}}).
\end{eqnarray}
In essence, we assume that the indivisibility of the mixture
components is not too small compared to the clustering difficulty.

With this notation, the following result applies to i.i.d. labeled
samples $\{(X_i,{Z}_i)\}_{i=1}^{n}$ from a
$K$-component mixture $\bar{\mathbb{P}}$.
%

%
\begin{theorem}[(Finite-sample angular structure)]
\label{ThmMain}
There are numbers ${c}, {c}_0,\break  {c}_1,{c}_2$
depending only on $b$ and $r$ such that for any
\mbox{$\delta\in(0,\frac{\llVert k \rrVert_{{\bar
{\mathbb{P}}}}}{b\sqrt{2\pi}})$}
satisfying condition~(\ref{EqnDevCondition}) and any
$t> \frac{{{c}_0}}{w_{\min}^3} \sqrt
{\varphi_n(\delta)}$, the
embedded data set $\{\Phi_{\mathcal V}(X_i),{Z}_i\}
_{i=1}^{n}$
has \mbox{$(\alpha,\theta)$-OCS} with
%
%
\begin{equation}
\label{EqnFiniteSampleGuarantee} \llvert\cos\theta\rrvert\leq\frac
{{{c}_0}\sqrt{\varphi_n
(\delta)}}{w_{\min}^3
t- {{c}_0}\sqrt{\varphi_n(\delta)}} \quad\mbox
{and}\quad\alpha\leq\frac{{c}_1}{(w_{\min})^{3/2}} \varphi
_n(\delta) + \psi(2t),\hspace*{-20pt}
\end{equation}
and this event holds probability at least $1 -
8K^2\exp({\frac{-{c}_2 n\delta^4}{ \delta^2 +
\mathcal{S}_{\max}({\bar{\mathbb{P}}}) + {{\mathcal{C}}({{\bar
{\mathbb{P}}}})} }} )$.
\end{theorem}
%

%
\begin{figure}[t]

\includegraphics{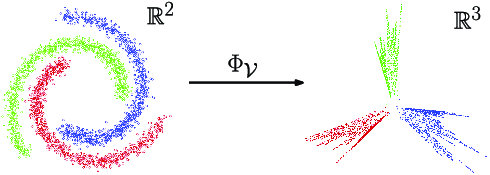}

\caption{According to Theorem~\protect\ref{ThmMain}, the normalized
Laplacian embedding of
i.i.d. samples from a nonparametric mixture with small overlap,
indivisible components and large enough sample size, has
$(\alpha,\theta)$-OCS with $\alpha\ll1$
and $\theta\ll1$. The left plot shows i.i.d. samples in
$\mathbb{R}^2$, and the right plot displays the image (in $\mathbb
{R}^3$) of
these data under the normalized Laplacian embedding, $\Phi_{\mathcal
V}$. The
embedding was performed using a regularized Gaussian kernel. The color
of each point indicates the latent label of
that point.}\label{FigThm}
\end{figure}

Theorem~\ref{ThmMain} establishes that the embedding of
i.i.d. samples from a finite mixture ${\bar{\mathbb{P}}}$ has
orthogonal cone
structure (OCS) \textit{if} the components have small overlap and good
indivisibility. This result holds with high probability on the
sampling from~${\bar{\mathbb{P}}}$. See Figure~\ref{FigThm} for an
illustration of
the theorem.


%
\begin{figure}[b]

\includegraphics{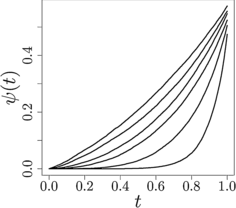}

\caption{The tail decay function $\psi(t)$ roughly
follows a power law for the standard Gaussian distribution and Gaussian
kernel with bandwidths $\nu\in\{0.15, 0.45, 0.75, 1, 1.5,2.5\}$.}
\label{FigTail}
\end{figure}

The tail decay of the mixture components enters the
bounds on $\alpha$ and $\theta$ in different ways: the bound
on $\theta$ is inversely proportional to
$t$, but the bound on $\alpha$ is tighter for smaller
$t$.
Depending on how quickly $\psi$ increases with $t$, it
may very well be the dominant term in the bound on $\alpha$.
For example, if there is a $\gamma>0$ such that $\psi(t)
\le t^\gamma$ for all $t\in(0,1)$, and we set
$t= \varphi_n^\beta$ for some $\beta\in
(0,1)$, then we obtain the simplified bounds
\[
\llvert\cos\theta\rrvert\leq\frac{c}{w_{\min}^3} \varphi_n
^{1/2 -\beta}\quad\mbox{and}\quad\alpha\leq\frac{c}{(w_{\min})^{3/2}} \bigl(
\varphi_n+ \varphi_n^{\gamma\beta} \bigr).
\]
Indeed, we find that whenever $\gamma\beta< 2$, the tail decay
function is the dominant term in the bound on~$\alpha$. Note that
this power law tail decay is easy to verify for the Gaussian
distribution with Gaussian kernel from Example~\ref{ExGaussian}; see
Figure~\ref{FigTail}.


Finally, the numbers $c,c_0,c_1,c_2$ increase as
the kernel bound $b$ increases and as $r$ decreases.
This is where we need the tail truncation condition $r>0$.
This assumption is common in the literature; see Cao and Chen
\cite{compharmonic}, for example. Both von Luxburg, Belkin and Bousquet
\cite{luxburg} and Rosasco, Belkin and De Vito \cite{rosasco} assume
$k(x,y) \ge r>0$, which is more restrictive. Note that
this automatically holds if we add a positive constant to any kernel.
This is sometimes called \textit{regularization} and can significantly
increase the performance of spectral clustering in practice
\cite{chen}.


\subsection{Algorithmic consequences}\label{SecAlg}

In this section we apply our theory to study the performance of
spectral clustering. The standard spectral clustering algorithm
applies $K$-means to the embedded dataset. For completeness, we
give pseudo code for the update step of $K$-means in Algorithm \ref{AlgAngularKMeans} below.

\begin{algorithm}[b]
\caption{$K$-means update}\label{AlgAngularKMeans}
\begin{algorithmic}[h!]
\State{\bf Input:} Normalized embedded data $y_i: =
\frac{\Phi_{\mathcal V}(X_i)}{\llVert\Phi_{\mathcal V}(X_i)
\rrVert}$ for $i
= 1,\ldots, {n}$, and mean vectors
$\{\mathbf a_1,\ldots,\mathbf a_K\}$
\For{$m\in\{1,\ldots, K\}$} \State
\[
{\hat{\mathcal{Z}}}_m \leftarrow\Bigl\{i\dvtx  m= \operatorname{argmin}\limits_{\ell}
\llVert{\mathbf a_\ell- y_i} \rrVert\Bigr\}
\]
\State
\[
\mathbf a'_m\leftarrow\sum
_{i\in
{\hat{\mathcal{Z}}}_m} \frac{y_i} {\llvert{\hat{\mathcal
{Z}}}_m\rrvert}
\]
\EndFor
\State\Return$\{{\hat{\mathcal{Z}}}_1, \ldots, {\hat{\mathcal
{Z}}}_K\}$ and
$\{{\mathbf a'_1},\ldots,{\mathbf a'_K}\}$
\end{algorithmic}
\end{algorithm}

In practice, we have found that applying $K$-means to an embedded
dataset works well if the underlying orthogonal cone structure is
``nice enough.'' The following proposition provides a quantitative
characterization of this phenomenon. It applies to an embedded data
set $\{\Phi_{\mathcal V}(X_i),{Z}_i\}
_{i=1}^{n}$ with
$(\alpha,\theta)$-OCS, and an initialization of
$\mathbf a_1,\ldots,\mathbf a_K$ as uniformly random
orthonormal vectors. Recall the notation ${\mathcal{Z}}_m= \{i \in
[{n}] |Z_i = m\}$.

%
\begin{proposition}\label{PropSpike}
Suppose $\theta$ and $\alpha$ are sufficiently small that
%
%
\begin{eqnarray}\label{AssIneq}
\frac{ \alpha{n}+ (1-\alpha) \llvert {\mathcal{Z}}_m
\rrvert \sin
\theta}{(1-\alpha) \llvert {\mathcal{Z}}_m \rrvert }&\le&\sin\frac{\pi
}{8}\quad\mbox{and}
\nonumber\\[-8pt]\\[-8pt]\nonumber
\frac{(1-\alpha)\llvert {\mathcal{Z}}_m \rrvert\cos\theta- \alpha
n}{\llvert {\mathcal{Z}}_m \rrvert + \alpha n} &\ge& \frac{1} 2, \qquad
m= 1,\ldots,K.
\end{eqnarray}
Then there is a constant $c_K$ such that with probability at least $1
- \frac{4 c_K \theta}{2 \pi}$ over the random initialization, the
$K$-means algorithm misclusters at most $\alpha n$ points.
When $K = 2$, we have $c_K = 1$.
\end{proposition}

Intuitively, condition~(\ref{AssIneq}) requires $\alpha$ and
$\theta$ to be small enough so that the different cones from the
$(\alpha,\theta)$-OCS do not overlap.

\begin{pf*}{Proof of Proposition \ref{PropSpike}}
We provide a detailed proof for the case $K= 2$. By the
definition of $(\theta,\alpha)$-OCS, there exist orthogonal
vectors $e_1,e_2$ such that a fraction $1-\alpha$ of the embedded
samples with latent label $m$ lie within an angle $\theta$ of
$e_m$, $m= 1,2$. Let us say that the initialization is
\emph{unfortunate} if some $\mathbf a_j$ falls within angle
$\frac{\theta}{2}$ of the angular bisector of $e_1$ and $e_2$, an
event which occurs with probability $\frac{4\theta}{2\pi}$.

Suppose without loss of generality that $\mathbf a_1$ is closer to
$e_1$, and let $\mathbf a'_1,\mathbf a'_2$ denote the updates
\[
\mathbf a'_m= \sum_{i\in{\hat{\mathcal{Z}}}_m}
\frac
{v_i} {\llvert{\hat{\mathcal{Z}}}_m\rrvert},\qquad m = 1,2.
\]
If the initialization is \emph{not} unfortunate, then all points in
the $\theta$-cone around $e_1$ are closer to $\mathbf a_1$ than
$\mathbf a_2$. In this case, the $(\theta,\alpha)$-OCS
implies that the $e_2$-coordinate of~$\mathbf a'_1$ is at most
\[
\frac{ \alpha{n}+ (1-\alpha) \llvert {\mathcal{Z}}_1 \rrvert
\sin\theta}{(1-\alpha) \llvert {\mathcal{Z}}_1 \rrvert } \le\sin\frac
{\pi}{8},
\]
and the $e_1$-coordinate of $\mathbf a'_1$ is at least
\[
\frac{(1-\alpha)\llvert {\mathcal{Z}}_1 \rrvert\cos\theta- \alpha
n}{\llvert {\mathcal{Z}}_1 \rrvert + \alpha n} \ge\frac{1} 2.
\]
We conclude that all points in the $\theta$-cone about $e_1$ are
closer to ${\mathbf a'_1}$ than ${\mathbf a'_2}$. Consequently,
we find that after a single update step of $K$-means, all but a
fraction $\alpha$ of the samples are correctly labeled. Moreover,
this holds for all subsequent $K$-means updates. This completes
the proof for $K= 2$.

The proof for general $K$ follows the same steps. The probability
that any $\mathbf a_m$ falls within angle $\frac{\theta}{2}$ of
the angular bisector of any pair $e_j,e_\ell$ is still proportional to~$\theta$, with a constant of proportionality $c_K$ that depends on
$K$.
\end{pf*}


\section{Proofs}
\label{SecProofs}

We now turn to the proofs of our main results, beginning with the
population level result stated in Theorem~\ref{ThmPop}. We then
provide the proof of Theorem~\ref{ThmMain} and
Proposition~\ref{PropSpike}.


\subsection{Proof of Theorem~\texorpdfstring{\protect\ref{ThmPop}}{1}}

Our proof leverages an operator perturbation theorem due to
Stewart~\cite{stewart} to show that ${\mathcal{Q}}$ is an approximate
invariant subspace of the normalized Laplacian operator
$\bar{\mathbf{T}}$ from equation~(\ref{EqnDefnNormLapOp}). Recalling
that ${\Pi_{{\mathcal{Q}}}}$ denotes the projection onto subspace
${\mathcal{Q}}$ (with
${\Pi_{{{\mathcal{Q}}^\perp}}}$ defined analogously), consider the
following three
operators:
\[
{\mathbf{A}}: ={{\Pi_{{\mathcal{Q}}}} \bar{\mathbf{T}} {\Pi
^*_{{\mathcal{Q}}}}},
\qquad{\mathbf{B}}: = {{\Pi_{{{\mathcal{Q}}^\perp}}} \bar{\mathbf{T}}
{\Pi
^*_{{{\mathcal{Q}}^\perp}}}}\quad\mbox{and}\quad{\mathbf{G}}: = {{\Pi
_{{{\mathcal{Q}}^\perp}}}
\bar{\mathbf{T}} {\Pi^*_{{\mathcal{Q}}}}}.
\]
By definition, a subspace ${\mathcal{Q}}$ is invariant under $\bar
{\mathbf{T}}$
if and only if ${\mathbf{G}}= 0$. In our setting, this ideal situation
occurs when there is no overlap between mixture components. More
generally, operator perturbation theory can be used to guarantee that a
space is approximately invariant as long as the Hilbert--Schmidt norm
${\vvvert {\mathbf{G}} \vvvert_{\mathrm{HS}}}$ is not too large relative to the spectral separation
between ${\mathbf{A}}$ and ${\mathbf{B}}$. In particular, define the
quantities
\[
\gamma: =\vvvert{\mathbf{G}} \vvvert_{\mathrm{HS}}\quad\mbox{and}\quad{
\operatorname{sep}({\mathbf{A}}, {\mathbf{B}})}: =\inf\bigl\{
\llvert a - b
\rrvert|a \in\sigma({\mathbf{A}}), b \in\sigma({\mathbf{B}}) \bigr\}.
\]
In application to our problem, Theorem 3.6 of Stewart~\cite{stewart}
guarantees that as long as $\frac{\gamma}{{\operatorname
{sep}({\mathbf{A}}, {\mathbf{B}})}} <
\frac{1}{2}$, then there is an operator ${\mathbf{S}}\dvtx  {\mathcal
{Q}}\rightarrow
{{\mathcal{Q}}^\perp}$ such that
%
%
\begin{eqnarray}
\label{EqnMyStewartBound} \vvvert{\mathbf{S}} \vvvert
_{\mathrm{HS}} & \leq&
\frac{2 \gamma}{{\operatorname{sep}({\mathbf
{A}}, {\mathbf{B}})}}
\end{eqnarray}
such that ${\operatorname{Range}}({\Pi^*_{{\mathcal{Q}}}} + {\Pi
^*_{{{\mathcal{Q}}^\perp}}} {\mathbf{S}}
)$ is an invariant subspace of $\bar{\mathbf{T}}$.

Accordingly, in order to apply this result, we first need to control
the quantities $\vvvert {\mathbf{G}} \vvvert_{\mathrm{HS}}$ and ${\operatorname{sep}({\mathbf{A}},
{\mathbf{B}})}$. The bulk of
our technical effort is devoted to proving the following two lemmas:
%

\begin{lemma}[(Hilbert--Schmidt bound)]\label{LemGammaBound}
We have
%
%
\begin{eqnarray}
{\vvvert{\mathbf{G}} \vvvert_{\mathrm{HS}}} &\leq& \frac{\sqrt{K(12+ b_{\max})}}{w_{\min}} \sqrt{
\mathcal{S}_{\max}({\bar{\mathbb{P}}}) + {{\mathcal{C}}({\bar{
\mathbb{P}}})}},
\end{eqnarray}
where $b_{\max}: =\max_{m= 1, \ldots, K}
{\llVert\int k_m(x,y)\,d\mathbb{P}_m(y) \rrVert_\infty^2}$.
\end{lemma}


%
%
\begin{lemma}[(Spectral separation bound)]\label{LemSpecSep}
Under the hypothesis of the theorem, we have
\begin{subequations}
\begin{eqnarray*}
\sigma_{\min}({\mathbf{A}}) & \geq& 1
- 13 K\bigl[\mathcal{S}_{\max}({
\bar{\mathbb{P}}}) + {{\mathcal{C}}({\bar{\mathbb{P}}})}\bigr]^{1/2}
\end{eqnarray*}
and
\begin{eqnarray*}
\sigma_{\max}({\mathbf{B}}) & \leq&1 - \frac{\Gamma^2} 8 +
\frac{3
[{\mathcal{S}_{\max}({\bar{\mathbb{P}}}) + {{\mathcal
{C}}({\bar{\mathbb{P}}})}}]^{1/2}}{w_{\min}}.
\end{eqnarray*}
\end{subequations}
Consequently, the spectral separation is lower bounded as
%
%
\begin{eqnarray}
{\operatorname{sep}({\mathbf{A}}, {\mathbf{B}})} & \geq&\frac
{\Gamma^2}{16}.
\end{eqnarray}
\end{lemma}

See Appendices~\textup{A.1}~and~\textup{A.2} (supplementary material \cite{SchiebingerWainwrightYu14Sup}), respectively, for the proof of these lemmas.

Combined with our earlier bound~(\ref{EqnMyStewartBound}),
these two lemmas guarantee that
%
%
\begin{eqnarray}
\label{EqnExplicitStewart} {\vvvert{\mathbf{S}} \vvvert
_{\mathrm{HS}}} & \leq&16\sqrt{12+
{b_{\max}}} \varphi({\bar{\mathbb{P}}},k).
\end{eqnarray}

Moreover,\vspace*{-1pt} we find that ${\operatorname{Range}}({\Pi^*_{{\mathcal
{Q}}}} + {\Pi^*_{{{\mathcal{Q}}^\perp}}} {\mathbf{S}})$ is equal
to~${\mathcal{R}}$, the principal eigenspace
of $\bar{\mathbf{T}}$. Indeed, by Stewart's theorem, the spectrum of
$\bar{\mathbf{T}}$
is the disjoint union
$\sigma(\bar{\mathbf{T}}) = \sigma({\mathbf{A}}+ {\mathbf
{G}}^*{\mathbf{S}}) \cup\sigma({\mathbf{B}}-
{\mathbf{S}}{\mathbf{G}}^*)$.
After some calculation using the upper bound on $\varphi({{\bar
{\mathbb{P}}}},k)$ in the theorem hypothesis, we find that the spectrum of $\bar
{\mathbf{T}}
$ satisfies $\sigma_{\min}({\mathbf{A}}+ {\mathbf{G}}^* {\mathbf
{S}}) > \sigma_{\max
}({\mathbf{B}}- {\mathbf{S}}{\mathbf{G}}^*)$,
and any element $x \in{\operatorname{Range}}({\Pi^*_{{\mathcal
{Q}}}} + {\Pi^*_{{{\mathcal{Q}}^\perp}}} {\mathbf{S}})$ must satisfy
\[
\sup_{q\in{\mathcal Q}} \biggl\{\frac{x^* \bar{\mathbf{T}}x}{x^* x} \bigg|x
= \bigl({
\Pi^*_{{\mathcal{Q}}}} + {\Pi^*_{{{\mathcal{Q}}^\perp}}} {\mathbf
{S}}\bigr)q \biggr\} >
\sigma_{\max}\bigl({{\mathbf{B}}- {\mathbf{S}} {\mathbf{G}}^*}\bigr).
\]
Therefore ${\operatorname{Range}}({\Pi^*_{{\mathcal{Q}}}} + {\Pi
^*_{{{\mathcal{Q}}^\perp}}} {\mathbf{S}}
) = {\mathcal{R}}$.

The only remaining step is to translate
bound~(\ref{EqnExplicitStewart}) into a bound on the norm
$\vvvert {\Pi_{{\mathcal{R}}}} - {\Pi_{{\mathcal{Q}}}} \vvvert
_{\mathrm{HS}}$. Observe that the
difference
of projection operators can be written as
\[
{\Pi_{\mathcal{R}}} - {\Pi_{{\mathcal Q}}} = ({\Pi_{\mathcal{R}}} + {
\Pi_{\mathcal{R}^\perp}}) ({\Pi_{\mathcal{R}}} - {\Pi_{{\mathcal
Q}}}) = {\Pi_{\mathcal{R}}} {\Pi_{{\mathcal Q}^\perp}} - {\Pi_{\mathcal
{R}^\perp}} {
\Pi_{{\mathcal Q}}}.
\]
Now Lemma 3.2 of Stewart~\cite{stewart} gives the explicit representations
\[
\Pi_{{\mathcal{R}}} = \bigl({\mathbf{I}}+ {\mathbf{S}}^*{\mathbf{S}}
\bigr)^{-1/2} \bigl({\Pi_{{\mathcal{Q}}}} + {\mathbf{S}}^* {
\Pi_{{{\mathcal{Q}}^\perp}}}\bigr)
\]
and
\[
\Pi_{{\mathcal
{R}}^\perp} = \bigl({
\mathbf{I}}+ {\mathbf{S}} {\mathbf{S}}^*\bigr)^{-1/2} ({
\Pi_{{{\mathcal{Q}}^\perp}}} + {\Pi_{{\mathcal{Q}}}} {\mathbf{S}}).
\]
%
Consequently, we
have
\[
\vvvert{\Pi_{\mathcal{R}^\perp}} {\Pi_{{\mathcal{Q}}}} \vvvert
_{\mathrm{HS}} \le
\vvvertdu\bigl({\mathbf{I}}+ {\mathbf{S}} {\mathbf{S}}^*\bigr)^{-1/2} {\mathbf{S}} \vvvertdu_{\mathrm{HS}}
\]
and
\[
\vvvert{\Pi_{\mathcal{R}}} {\Pi_{{{\mathcal{Q}}^\perp}}} \vvvert
_{\mathrm{HS}} \le \vvvertdu
\bigl({\mathbf{I}}+ {\mathbf{S}}^* {\mathbf{S}}\bigr)^{-1/2} {
\mathbf{S}}^* \vvvertdu_{\mathrm{HS}}.
\]
By the continuous functional calculus (see Section~VII.1
of Reed and Simon~\cite{reedsimon}), we have the expansion
\begin{eqnarray*}
\bigl({\mathbf{I}}+ {\mathbf{S}} {\mathbf{S}}^*\bigr)^{-1/2} & =& \sum
_{n=1}^\infty\pmatrix{2n
\cr
n}
\frac{({\mathbf{S}}{\mathbf{S}}^*)^{n-1}}{2^{2n}}.
\end{eqnarray*}
Putting together the pieces, in terms of the shorthand ${\varepsilon}=
{\vvvert {\mathbf{S}} \vvvert_{\mathrm{HS}}}$, we have
\begin{eqnarray*}
\vvvert{\Pi_{{\mathcal{R}}}} - {\Pi_{{\mathcal{Q}}}} \vvvert
_{\mathrm{HS}} & \leq&\frac{{\varepsilon}}{2} \sum_{n=1}^\infty
\pmatrix{2n
\cr
n} \biggl( \frac{{\varepsilon}}{2} \biggr)^{2(n-1)} =
\frac{2}{{\varepsilon}} \biggl(\frac{1}{\sqrt{1-{\varepsilon}^2}} - 1 \biggr)
\leq{\varepsilon},
\end{eqnarray*}
which completes the proof.


\subsection{Proof of Theorem~\texorpdfstring{\protect\ref{ThmMain}}{2}}

We say that a $K$-element subset (or \mbox{$K$-tuple}) of
$\{X_1,\ldots,X_{n}\}$ is \textit{diverse} if the
latent labels
of all points in the subset are distinct. Given some $\theta\in
(0,\frac{\pi}{4})$, a $K$-tuple is $\theta$-orthogonal if all
its distinct pairs, when embedded, are orthogonal up to angle
$\frac{\theta}{2}$. In order to establish
$(\alpha,\theta)$-angular structure in the normalized Laplacian
embedding of $\{X_1,\ldots,X_{n}\}$, we must show that
there is a subset of $\{X_1,\ldots,X_{n}\}$ with at
least $(1-\alpha){n}$ elements, and with the property that
every diverse $K$-tuple from the subset is $\theta$-orthogonal.

We break the proof into two steps. We first lower bound the total
number of $K$-tuples that are diverse and $\theta$-orthogonal.
In the second step we construct the desired subset. We present the
first step below and defer the second step to
Appendix~\textup{B.1} in the supplementary material \cite{SchiebingerWainwrightYu14Sup}.

\subsubsection*{Step 1}
Consider a diverse $K$-tuple $(X^1,\ldots,X^K)$
constructed randomly by selecting $X^m$ uniformly at
random from the set $\{X_i|{Z}_i= m\}$
for $m= 1,\ldots,K$. Form the $K\times K$ random
matrix
\begin{eqnarray*}
{V}&=& \lleft[\matrix{ |&&|
\cr
\Phi_{\mathcal V}\bigl(X^1
\bigr) & \cdots& \Phi_{\mathcal V}\bigl(X^K\bigr)
\cr
| & & | }
\rright],
\end{eqnarray*}
where $\Phi_{\mathcal V}$ denotes the normalized Laplacian embedding from
equation~(\ref{EqnSpectralEmbedding}). Let ${\tilde{V}}$ denote an
independent copy of ${V}$. Let ${Q}\in\mathbb{R}^{K\times
K}$
denote the diagonal matrix with entries ${Q}_{mm} =
\frac{q_m(X^m)}{\llVert q_m \rrVert_{{{\bar{\mathbb
{P}}}_n}}}$, where ${\bar{\mathbb{P}}_n}$ is the
empirical distribution over the samples
$X_1,\ldots,X_{n}$, and define $Q_{\max}: =
\max_m\frac{{\llVert q_m \rrVert_\infty}}{\llVert q_m
\rrVert_{\bar{\mathbb{P}}}}$.\vspace*{2pt}

At the\vspace*{1pt} core of our proof lies the following claim involving a constant
$c_3$. For at least a fraction $1-\frac{2Kc_3
\varphi_n(\delta)} {\sqrt{w_{\min}}}$ of the diverse
$K$-tuples, we have the inequality
%
%
\begin{equation}
\label{EqnKeyClaim} \vvvertdu{V^T} {\tilde{V}}- {Q}^2
\vvvertdu_{\mathrm{HS}} \leq\frac{32
\sqrt{3}}{w_{\min}^3} Q_{\max}^2
\sqrt{c_3}\sqrt{\varphi_n(\delta)},
\end{equation}
holding on a high probability set ${\mathcal A}$. For the moment, we take
this claim as given, before returning to define ${\mathcal A}$
explicitly and
prove the claim.

When inequality~(\ref{EqnKeyClaim}) is satisfied, we obtain the
following upper bound on the off-diagonal elements of ${{{V}}^T}
{\tilde{V}}$:
\[
\bigl({{{V}}^T} {\tilde{V}} \bigr)_{m\ell}\le
\frac{32 \sqrt{3} Q_{\max}^2 }{w_{\min}^3} \sqrt{c_3 }\sqrt{\varphi
_n(\delta)}
\qquad\mbox{for }m\neq\ell.
\]
This is
useful because
\[
\cos \operatorname{angle}\bigl( \Phi_{\mathcal V}\bigl(X^m\bigr),
\Phi_{\mathcal V}\bigl(X ^\ell\bigr)\bigr) = \frac{({{{V}}^T}{\tilde
{V}})_{m\ell}}{\sqrt{({{{V}}^T}
{V})_{mm} ({{{V}}^T}{\tilde{V}})_{\ell\ell}}}.
\]
However, we must also lower bound
$\min_m({{{V}}^T}{\tilde{V}})_{mm}$.
To this end, by union bound, we obtain
%
%
\begin{eqnarray}\label{EqnTailBnd}
\Pr  \Bigl\{ \min
_m{Q}_{mm}^2 \le t \Bigr\} & =&
\Pr  \biggl\{ \min_m\frac{q_m^2(X^m)}{{\llVert q_m
\rrVert^2_{{\bar{\mathbb{P}}_n}}}} \le t
\biggr\} \nonumber
\\
&\le& \sum_{m= 1}^K\Pr
\biggl\{ \frac{q_m^2(X ^m)}{{\llVert q_m
\rrVert^2_{{\bar{\mathbb{P}}_n}}}} \le t \biggr\}
\\
&: =&\psi_n(t).\nonumber
%
\end{eqnarray}
On the set ${\mathcal A}_\psi: =\{ \sup_t \llvert \psi_n(t) -
\psi(t) \rrvert \le\delta\} \subset{\mathcal A}$, we may combine
equations~(\ref{EqnKeyClaim}) and~(\ref{EqnTailBnd}) to obtain
\[
\min_m \bigl({{{V}}^T} {\tilde{V}}
\bigr)_{mm}\ge t- \frac{32 \sqrt{3} Q_{\max}^2}{w_{\min}^3} \sqrt{c_3
}\sqrt{
\varphi_n(\delta)},
\]
with probability at least \mbox{$1 - \psi(2t)$.}
Therefore, there is a
$\theta$ satisfying
\[
\llvert\cos\theta\rrvert\le\frac{
32 \sqrt{3} Q_{\max}^2\sqrt{c_3}\sqrt{\varphi_n
(\delta)}}{w_{\min}^3t-
32 \sqrt{3} Q_{\max}^2\sqrt{c_3}\sqrt{\varphi_n
(\delta)}}
\]
such that at least a fraction $1-\frac{2Kc_3
\varphi_n(\delta)} {
\sqrt{w_{\min}}} - \psi(2t)$ of the diverse
$K$-tuples are $\theta$-orthogonal on the set ${\mathcal A}$.
This establishes the finite sample bound~(\ref
{EqnFiniteSampleGuarantee}) with $c_0: =2 c_3$ and
$c_1: =32 \sqrt{3} Q_{\max}^2 \sqrt{c_3}$.

It remains to prove the intermediate claim~(\ref{EqnKeyClaim}).
Define the matrix
%
%
\begin{eqnarray}
A&: =& \lleft[\matrix{ \displaystyle\biggl\langle\frac{q_1}{\llVert q_1 \rrVert
_{{{\bar{\mathbb{P}}}_n}}},v_1
\biggr\rangle_{{\bar{\mathbb{P}}_n}} & \cdots& \displaystyle\biggl\langle\frac
{q_1}{\llVert q_1 \rrVert
_{{\bar{\mathbb{P}}_n}}},v_K
\biggr\rangle_{{{\bar{\mathbb
{P}}}_n}}
\cr
\vdots& \ddots& \vdots
\cr
\displaystyle\biggl\langle\frac{q_K}{\llVert q_K \rrVert_{{{\bar
{\mathbb{P}}}_n}}},v_1 \biggr\rangle_{{\bar{\mathbb{P}}_n}} &
\cdots& \displaystyle\biggl\langle\frac{q_K}{\llVert q_K \rrVert_{{{\bar
{\mathbb{P}}}_n}}},v_K \biggr
\rangle_{{\bar{\mathbb{P}}_n}} } \rright].
\end{eqnarray}
Note that the entries of $A{{A}^T}$ are
\[
\bigl(A{{A}^T} \bigr)_{m\ell} = \frac{ \langle
{\Pi_{{\mathcal V}}}q_m,{\Pi_{{\mathcal V}}}q_\ell
\rangle_{{\bar{\mathbb{P}}_n}}}{\llVert q_m \rrVert_{{\bar{\mathbb
{P}}_n}}\llVert q_\ell\rrVert
_{{\bar{\mathbb{P}}_n}}}.
\]
The off-diagonal elements satisfy
\[
\bigl(A{{A}^T} \bigr)_{m\ell} \le3(\hat\varphi+ \sqrt{\hat
\mathcal{S}_{\max}})\qquad\mbox{for } m\ne\ell,
\]
where $\hat\varphi= \max_m\frac{\llVert q_m- {\Pi_{{\mathcal
V}}} q_m \rrVert_{{\bar{\mathbb{P}}_n}}}{ \llVert q_m
\rrVert_{{\bar{\mathbb{P}}_n}}}$, and \mbox{$\hat{\mathcal
{S}}_{\max}= \max_{m\neq
\ell} \frac{{\llVert q_\ell\rrVert^2_{\mathbb
{P}_m^n}}}{{\llVert q_m \rrVert^2_{{{\bar{\mathbb
{P}}}_n}}}}$} (and\vspace*{1pt} $\mathbb{P}_m^n$ denotes
the empirical distribution for the samples with latent label $m$).
Similarly, the diagonal elements satisfy
$\llvert (A{{A}^T} )_{mm} -1 \rrvert\le
3\hat\varphi$. Putting together the pieces yields
${\vvvert A{{A}^T} - I \vvvert_{\mathrm{HS}}}^2 \le3K^2 (\hat\varphi+
\hat\mathcal{S}_{\max})$, which in turn implies
\[
\label{EqAbndUseful} \vvvertdu{\bigl(A{{A}^T}\bigr)}^{-1} - I
\vvvertdu_{\mathrm{HS}}^2 \le\frac
{3K^2
(\hat\varphi+ \sqrt{\hat\mathcal{S}_{\max}})}{1-3K^2 (\hat
\varphi+ \sqrt{\hat\mathcal{S}_{\max}})}.
\]

We now transform this inequality into one involving ${{{V} }^T}{\tilde{V}}$.
Write $B= A{V}$ and ${\tilde B}= A{\tilde{V}}$,
and note that
${{{V}}^T}{\tilde{V}}=
{{B}^T}{(A{{A}^T})}^{-1}{\tilde B}$.
Therefore, we
find that
\begin{eqnarray*}\label{EqnAlmostDone}
\vvvertdu{V^T} {
\tilde{V}}- {Q}^2 \vvvertdu_{\mathrm{HS}} & \le& \vvvertdu
{{B}^T} {\tilde B}- {Q}^2 \vvvertdu_{\mathrm{HS}} +
\vvvertdu{{B}^T} \bigl[{\bigl(A{{A}^T}\bigr)}^{-1}
-I \bigr]{\tilde B} \vvvertdu_{\mathrm{HS}}
\\
& \le& 3 \vvvert{Q} \vvvert_{\mathrm{HS}} \vvvert B-{Q} \vvvert
_{\mathrm{HS}} + {\vvvert B \vvvert_{\mathrm{HS}}}^2 \vvvertdu{
\bigl(A{{A}^T}\bigr)}^{-1} - I \vvvertdu_{\mathrm{HS}},
\end{eqnarray*}
where\vspace*{1pt} the last inequality used $\vvvert B \vvvert
_{\mathrm{HS}} \le2 \vvvert {Q} \vvvert
_{\mathrm{HS}}$.
Now note that the entries of $B$ are \mbox{$B_{m\ell
} =\frac{{\Pi_{{\mathcal V}}}q_m(X^\ell)}{\llVert q_m \rrVert_{{\bar
{\mathbb{P}}_n}}}$.}
Therefore the difference $B- {Q}$ satisfies
%
%
\begin{equation}
\label{EqnExpectationBS} \mathbb{E} \bigl[\vvvert B- {Q} \vvvert
_{\mathrm{HS}}^2
| X_1,\ldots,X_{n} \bigr] \le K^2 \biggl(
\frac{\hat\varphi}{\sqrt{{{\hat w}_{\min}}}} + \sqrt{\hat\mathcal
{S}_{\max} } \biggr)^2
+ K \frac{\hat\varphi^2}{{{\hat w}_{\min}}},
\end{equation}
where ${{\hat w}_{\min}}= \min_m\frac{{n}_m}{{n}}$,
and the
expectation above is over the selection of the random $K$-tuple
$(X^1,\ldots,X^K)$.\footnote{Note that there are two
different types of randomness at play in the construction of
${V}$ and hence~$B$; there is randomness in the generation of
the i.i.d. samples $X_1,\ldots,X_{n}$ from $\bar{\mathbb{P}}$, and
there is randomness in the selection of the diverse $K$-tuple
$(X^1,\ldots,X^K)$.}

Both $\hat\mathcal{S}_{\max}$ and $\hat\varphi$ are small with high
probability. Indeed, Bernstein's inequality guarantees that
%
%
\begin{equation}
\label{EqnSimilProb} \sqrt{\hat\mathcal{S}_{\max}} \le\sqrt{
\mathcal{S}_{\max}} + \delta
\end{equation}
with probability at least $1-2K^2 \exp
{\frac{-{n}(\mathcal{S}_{\max}+\delta^2)^2}{8Q_{\max}^2
(2\mathcal{S}_{\max}+
\delta^2)}}$.
We control $\hat\varphi$ with a finite sample version of
Theorem~\ref{ThmPop}, which we state as Proposition~\ref{PropFinite} below.

Let ${\mathcal V}= \operatorname{span}\{v_1,\ldots,v_K\}$ denote the
principal eigenspace of the normalized Laplacian matrix.

%
\begin{proposition}
\label{PropFinite}
There are constants ${c_2'}, c_3$ such that for any $\delta\in
(0, \frac{\llVert k \rrVert_{{\bar{\mathbb
{P}}}}}{b\sqrt{2\pi}})$
satisfying condition~(\ref{EqnDevCondition}), we have
%
%
\begin{equation}
\label{EqnFiniteOne} \hat\varphi\le c_3\varphi_n(\delta)
\end{equation}
with probability at least $1 - 10K\exp({ \frac{-{n}
{c_2'}\delta^4}{\delta^2+\mathcal{S}_{\max}({{\bar{\mathbb
{P}}}}) +
{{\mathcal{C}}({\bar{\mathbb{P}}})}}} )$.
\end{proposition}

See Section~\ref{SecPropFiniteProof} for the proof of this
auxiliary result.

On the set $\{\hat\varphi\le c_3\varphi_n(\delta)\}\cap\{
{{\hat w}_{\min}}
\ge\frac{1} 2 w_{\min}\}:={\mathcal A}_\zeta\cap{\mathcal A}_w$,
equation~(\ref{EqnExpectationBS}) simplifies to
\[
\mathbb{E} \bigl[\vvvert B- {Q} \vvvert_{\mathrm{HS}}^2 |
X_1,\ldots,X_{n} \bigr] \le\frac{4K^2
c_3^2\varphi_n^2(\delta)} { {w_{\min}}},
\]
whenever $ ( \frac{\sqrt2 c_3\varphi_n(\delta
)}{\sqrt{w_{\min}}}
+ \mathcal{S}_{\max}+ \delta)^2 \le
\frac{3c_3^2\varphi_n^2(\delta)}{w_{\min}}$, which
is a consequence
of condition~(\ref{EqnDevCondition}). By Markov's inequality we
obtain the following result: at least a fraction $1-\frac{2K
c_3\varphi_n(\delta)} { \sqrt{w_{\min}}}$ of the diverse
$K$-tuples satisfies
%
%
\begin{equation}
\label{EqnBS} {\vvvert B-{Q} \vvvert_{\mathrm{HS}}}^2\le
\frac{2Kc_3
\varphi_n} { \sqrt{w_{\min}}}.
\end{equation}
For the diverse $K$-tuples that do satisfy
inequality~(\ref{EqnBS}) we find that
\[
\vvvertdu{{{V}}^T} {\tilde{V}}- {Q}^2 \vvvertdu
_{\mathrm{HS}} \le\biggl( \frac{6
\sqrt{2}
K^{3/2} Q_{\max}} {w_{\min}^{1/4}} + 32 \sqrt{3} K^3
Q_{\max}^2 \biggr) \sqrt{c_3} {\sqrt{
\varphi_n}},
\]
valid on the set ${\mathcal A}= {\mathcal A}_{w} \cap{\mathcal A}_{q}
\cap{\mathcal A}_{\psi}
\cap\{\hat\varphi\le c_3\varphi_n(\delta)\} \cap\{ \sqrt
{\hat\mathcal{S}_{\max}} \le\sqrt{\mathcal{S}_{\max}} +
\delta\}$, thereby establishing the
bound~(\ref{EqnKeyClaim}).

To complete the first step of the proof of Theorem~\ref{ThmMain}, it
remains to control the probability of ${\mathcal A}$. By Hoeffding's
inequality, we have ${\mathbb{P}}[{\mathcal A}_w] \geq1 - K
e^{{-{n}w_{\min}^2}/2}$. Finally, an application of
Bernstein's inequality controls the probability of ${{\mathcal A}_q}$,
and an
application of Glivenko--Cantelli controls the probability of~${{\mathcal A}_\psi}$. Putting together the pieces we find that
${\mathcal A}$
holds with probability at least $ 1 - 8K^2
\exp({\frac{-{n}c_2\delta^4}{\delta^2 + {\mathcal{S}}_{\max
}({\bar{\mathbb{P}}}) +
{{\mathcal{C}}({\bar{\mathbb{P}}})}}} )$, where $c_2:= \min
({c_2'},
\frac{1} {8Q_{\max}^2} )$.


\subsection{Proof of Proposition~\texorpdfstring{\protect\ref{PropFinite}}{2}}\label{SecPropFiniteProof}

Consider the operator $\hat{\mathbf{T}}\dvtx  {L^2({{\bar{\mathbb
{P}}}_n})} \to{L^2({\bar{\mathbb{P}}_n})}$
defined by
\[
(\hat{\mathbf{T}}f) (x) = \int\frac{1}{\bar{q}_n(x)}k(x,y)\frac
{f(y)}{\bar{q}_n(y)} \,d{\bar{
\mathbb{P}}_n}(y),
\]
where $\bar{q}_n(x) = \frac{1}{{n}} \sum_{i= 1}^K
k(X_i,x)$ is the square root kernelized density for the
empirical distribution ${\bar{\mathbb{P}}_n}$ over the data
$X_1,\ldots,X_{n}$. 
${\bar k}^n(x,y):=\frac{k(x,y)}{\bar{q}_n(x)\bar{q}_n(y)}$ for the
normalized kernel
function. Note that for any $f \in{L^2({\bar{\mathbb{P}}_n})}$ and
$v\in\mathbb{R}^{n}$ with coordinates $v_i = f(X_i)$, we have
$(\hat{\mathbf{T}}f)(X_j) = (Lv)_j$, where $L$ is the
normalized Laplacian matrix~(\ref{EqnNormLapMat}). Consequently, the
principal eigenspace ${\mathcal V}$ of $L$ is isomorphic to the
principal eigenspace of $\hat{\mathbf{T}}$ which we also denote by
${\mathcal V}$ for simplicity.

To prove the proposition, we must relate $\hat{\mathbf{T}}$ to the
normalized Laplacian operator $\bar{\mathbf{T}}$. These operators
differ in
both their measures of integration---namely, ${\bar{\mathbb{P}}_n}$
versus $\bar{\mathbb{P}}$---and
their kernels, namely $\frac{k(x,y)}{\bar{q}_n(x)\bar{q}_n(y)}$ versus
$\frac{k(x,y)}{\bar{q}(x)\bar{q}(y)}$. To bridge the gap we
introduce an
intermediate operator $\tilde{\mathbf{T}}\dvtx {L^2({{\bar{\mathbb
{P}}}_n})}\to{L^2({\bar{\mathbb{P}}_n})}$
defined by
\[
(\tilde{\mathbf{T}}f) (x) = \int\frac{1}{\bar{q}(x)}k(x,y)\frac
{f(y)}{\bar{q}
(y)} \,d{
\bar{\mathbb{P}}_n}(y).
\]
Let ${\tilde{\mathcal V}}$ denote the principal eigenspace of
$\tilde{\mathbf{T}}$. The following lemma bounds the $\rho$-distance
between the principal eigenspaces of $\tilde{\mathbf{T}}$ and $\hat
{\mathbf{T}}$.

%
\begin{lemma}
\label{LemV}
For any $\delta\in[0, \frac{\llVert k \rrVert_{{\bar
{\mathbb{P}}}}}{b\sqrt{2\pi}}
]$
satisfying condition~(\ref{EqnDevCondition}), we have
%
%
\begin{equation}
\rho({\mathcal V},{\tilde{\mathcal V}}) \le\frac{{c_4}}{\Gamma^2}
\biggl(
\frac{1} {\sqrt{n}} +\delta\biggr),
\end{equation}
with probability at least $1- 4 e^{-n \pi \delta^2 / 4} - 2e^{-n \mathbf{E} \bar k \delta^2}$, where
${c_4}= 1024\sqrt{2\pi K}\frac{\llVert k \rrVert_{{\bar
{\mathbb{P}}}}b}{r^4}$.
\end{lemma}

See Appendix~\textup{B.3} (supplementary material \cite{SchiebingerWainwrightYu14Sup}) for a proof of this lemma.

We must upper bound $\llVert q_m- {\Pi_{{\mathcal V}}}q_m \rrVert
_{{\bar{\mathbb{P}}_n}}$. By the
triangle inequality,
\begin{eqnarray*}
\label{trianglebound}
\llVert q_m- {\Pi_{{\mathcal V}}}q_m
\rrVert_{{{\bar{\mathbb{P}}}_n}} &\le&\llVert q_m- {\Pi_{\mathcal{R}}}q_m
\rrVert_{{\bar{\mathbb{P}}_n}} + \llVert{\Pi_{\mathcal
{R}}}q_m- {
\Pi_{{\tilde{\mathcal V}}}}q_m \rrVert_{{{\bar
{\mathbb{P}}}_n}}
\\
&&{} +\llVert{\Pi_{{\mathcal V}}}q_m- {
\Pi_{{\tilde{\mathcal
V}}}}q_m \rrVert_{{\bar{\mathbb{P}}_n}}.
\end{eqnarray*}

Note that $\llVert{\Pi_{{\mathcal V}}}q_m- {\Pi_{{\tilde
{\mathcal V}}}}q_m \rrVert_{{\bar{\mathbb{P}}_n}}
\le\llVert q_m \rrVert_{{\bar{\mathbb{P}}_n}}\rho
({\tilde{\mathcal V}},{\mathcal V})$. We can control
this term with the lemma. The term $\llVert q_m- {\Pi_{\mathcal
{R}}}q_m \rrVert_{{\bar{\mathbb{P}}_n}}$ is the empirical
version of a quantity
controlled by Theorem~\ref{ThmPop}. We handle the empirical
fluctuations with a version of Bernstein's inequality. For $\delta
_{p}\ge
0$ we have the inequality
%
%
\begin{equation}
\label{EqnDev} \llVert q_m- {\Pi_{\mathcal{R}}}q_m
\rrVert_{{{\bar
{\mathbb{P}}}_n}} \le\llVert q_m- {\Pi_{\mathcal{R}}}
q_m \rrVert_{\bar{\mathbb{P}}} + \delta_{p}
\end{equation}
with probability at least $1-2\exp( -\frac{n
\delta_{p}^4}{8(\delta_{p}^2+{c_{\mathrm{pop}}}^2 \varphi^2) \tilde{Q}^2_{\max}} )$,
where
\[
\tilde{Q}_{\max}= \max_m\frac{{\llVert q_m- {\Pi
_{\mathcal{R}}}q_m \rrVert_\infty}}{\llVert q_m \rrVert_{\bar{\mathbb
{P}}}}\quad\mbox{and}\quad c_{\mathrm{pop}}: =16 \sqrt
{12 + \frac{b_{\max}}{K}}.
\]

It remains to control $\llVert{\Pi_{\mathcal{R}}}q_m- {\Pi
_{{\tilde{\mathcal V}}}}q_m \rrVert_{{\bar{\mathbb{P}}_n}}$.
Let ${\bar{\mathcal{H}}}$ denote the reproducing kernel Hilbert space
(RKHS)\footnote{We give a brief introduction to the theory of
reproducing kernel Hilbert spaces and provide some references for
further reading on the subject in Appendix~\textup{C.2} (supplementary material \cite{SchiebingerWainwrightYu14Sup}).} for the
kernel ${\bar k}$. Now we
define two integral operators on~${\bar{\mathcal{H}}}$. Let $\bar
{\mathbf{H}}$
denote the operator defined by
\[
(\bar{\mathbf{H}}h) (x) = \int{\bar k}(x,y) h(y) \,d\bar{\mathbb{P}}(y),
\]
and similarly let $\tilde{\mathbf{H}}\dvtx {\bar{\mathcal{H}}}\to{\bar
{\mathcal{H}}}$ denote the operator
defined by
\[
(\tilde{\mathbf{H}}h) (x) = \int{\bar k}(x,y) h(y) \,d{{\bar{\mathbb
{P}}}_n}(y).
\]
Both $\bar{\mathbf{H}}$ and $\tilde{\mathbf{H}}$ are self-adjoint,
compact operators
on ${\bar{\mathcal{H}}}$ and have real, discrete spectra.
Let ${\mathcal G}$ denote
the principal $K$-dimensional eigenspace of $\bar{\mathbf{H}}$, and let
${\tilde{\mathcal G}}$ denote the principal $K$-dimensional principal
eigenspace of
$\tilde{\mathbf{H}}$. The following lemma bounds the $\rho
$-distance between these subspaces of ${\bar{\mathcal{H}}}$.

%
\begin{lemma}
\label{LemG}
For any $\delta>0$ satisfying condition~(\ref{EqnDevCondition}), we have
\[
\rho({\mathcal G},{\tilde{\mathcal G}})\le\frac{c_5
}{\Gamma^2} \biggl(
\frac{1}{
\sqrt{n}} + \delta\biggr)
\]
with probability at least $1-2e^{-{n}\pi
\mathbb{E}{\bar k}({\bar X},{\bar X}) \delta^2 }$, where $c_5=
64\sqrt
{2\pi K}\sqrt{\mathbb{E}{\bar k}({\bar X},{\bar X})} \frac{b
}{r^2}$.
\end{lemma}

See Appendix~\textup{B.2} (supplementary material \cite{SchiebingerWainwrightYu14Sup}) for the proof of this
lemma.

By the triangle inequality, we have
\begin{eqnarray*}
\llVert{\Pi_{\mathcal{R}}}q_m- {\Pi_{{\tilde{\mathcal V}}}}q_m
\rrVert_{{\bar{\mathbb{P}}_n}} &\le&\llVert{\Pi_{\mathcal{R}}}q_m- {
\Pi_{{\mathcal G}}}q_m \rrVert_{{\bar{\mathbb{P}}_n}} + \llVert{
\Pi_{{\mathcal G}}}q_m- {\Pi_{{\tilde{\mathcal G}}}}q_m \rrVert
_{{\bar{\mathbb{P}}_n}}
\\
&&{} + \llVert{\Pi_{{\tilde{\mathcal G}}}}q_m- {
\Pi_{{\tilde{\mathcal
V}}}}q_m \rrVert_{{\bar{\mathbb{P}}_n}}.
\end{eqnarray*}
We claim that
%
%
\begin{equation}
\label{EqnZeroClaim} \llVert{\Pi_{\mathcal{R}}}q_m- {
\Pi_{{\mathcal G}}}q_m \rrVert_{{\bar{\mathbb{P}}_n}} = 0\quad\mbox{and}
\quad\llVert{\Pi_{{\tilde{\mathcal G}}}}q_m- {\Pi_{{\tilde{\mathcal V}}}}q_m
\rrVert_{{\bar{\mathbb{P}}_n}} = 0.
\end{equation}
We take these identities as given for the moment, before returning to
prove them at the end of this subsection.

Now the term $\llVert{\Pi_{{\mathcal G}}}q_m- {\Pi_{{\tilde
{\mathcal G}} }}q_m \rrVert_{{\bar{\mathbb{P}}_n}}$
can be controlled using the lemma in the following way. For any $h\in
{\bar{\mathcal{H}}}$, note that
\[
\llVert h \rrVert_{{\bar{\mathbb{P}}_n}}^2= \frac{1} n \sum
_{i= 1}^{n} \langle h,{\bar
k}_{X_i} \rangle_{{\bar{\mathcal
{H}}}}^2\le\frac{1} n
\sum_{i=
1}^{n} \llVert h \rrVert
_{\bar{\mathcal{H}}}^2 {\bar k}(X_i,X_i)
\]
by Cauchy--Schwarz for the RKHS inner product. Using this logic with
$h = {\Pi_{{\mathcal G}}}q_m- {\Pi_{{\tilde{\mathcal
G}}}}q_m$, we
find
%
%
\begin{equation}
\label{EqnVV} \llVert{\Pi_{{\mathcal G}}}q_m- {
\Pi_{{\tilde{\mathcal G}}}}q_m \rrVert_{{\bar{\mathbb{P}}_n}} \le
\llVert
q_m \rrVert_{\bar{\mathcal{H}}} \sqrt{\frac{1}{
n}\sum
_{i= 1}^{n} {\bar k}(X_i,X_i)}
\rho({\mathcal G},{\tilde{\mathcal G}}).
\end{equation}
Collecting our results and applying Lemmas~\ref{LemV} and~\ref{LemG}
yields
\[
\llVert q_m- {\Pi_{{\mathcal V}}}q_m \rrVert
_{{{\bar
{\mathbb{P}}}_n}} \le(c_{\mathrm{pop}}
\varphi+ \delta_{p} ){
\llVert q_m \rrVert_{\bar{\mathbb{P}}}} + \frac{c_n\llVert q_m
\rrVert_{{{\bar{\mathbb
{P}}}_n}}}{\Gamma^2} \biggl(
\frac{1}{
\sqrt
{n}} + \delta\biggr),
\]
where
%
%
\begin{eqnarray}
c_n&: =&\frac{256 \sqrt{2\pi K} b}{r^2} \biggl[ \frac{\llVert q_m
\rrVert_{\bar{\mathcal{H}}}}{\llVert
q_m \rrVert_{\bar{\mathbb{P}}}}\mathbb{E} {
\bar k}({\bar X},{\bar X}) + \frac{2b} {r^2} \biggr].
\end{eqnarray}
By an application of Bernstein's inequality, we have
\[
\llVert q_m \rrVert_{{\bar{\mathbb{P}}_n}} \le\sqrt{2}\llVert
q_m \rrVert_{\bar{\mathbb{P}}}
\]
with probability at least $1-2e^{{-{n}}/(16 Q_{\max}^2)}$.
For $\delta\in(0, \frac{1} {2 \sqrt{2 \pi} Q_{\max}})$, we have
\begin{eqnarray*}
&& 2 e^  {(- {n}{c_{\mathrm{pop}}}^2\delta^4)/(8\Gamma^4
\tilde{Q}_{\max}^2(\delta^2 + \mathcal{S}_{\max}({{\bar{\mathbb
{P}}}}) + {{\mathcal{C}}({\bar{\mathbb{P}}})}) )} + 6e^{-{n}\pi
\delta^2/2} + 2e^{-{n}/(16Q_{\max}^2)}
\\
&&\qquad  \leq 10 e^ {-({n
{c_2'}
\delta^4})/({\delta^2+ \mathcal{S}_{\max}({\bar{\mathbb{P}}}) +
{{\mathcal{C}}({\bar{\mathbb{P}}})}})},
\end{eqnarray*}
where $\delta_{p}= \frac{c_{\mathrm{pop}}\delta}{\Gamma^2}$, and
${c_2'}
= \min( \frac{c_{\mathrm{pop}}^2}{8\Gamma^4\tilde{Q}_{\max}^2},
\frac{\pi}2 )$. Modulo the claim, this proves the proposition
with $c_3= 2\max(c_{\mathrm{pop}},c_n)$.

We now return to prove claim~(\ref{EqnZeroClaim}). Note the
following relation between the eigenfunctions of $\bar{\mathbf{T}}$
and those of
$\bar{\mathbf{H}}$: if $r_i$ is an eigenfunction of $\bar{\mathbf
{T}}$ with
eigenvalue $\lambda_i$ and $\llVert r_i \rrVert_{{\bar
{\mathbb{P}}}}=1$, then
$g_i: ={\sqrt{\lambda_i}} r_i$ has unit
norm in ${\bar{\mathcal{H}}}$ and is an eigenfunction of $\bar
{\mathbf{H}}$ with eigenvalue
$\lambda_i$. Note\vspace*{1pt} that the eigenfunctions $r_i$ of
$\bar{\mathbf{T}}$ form an orthonormal basis of ${L^2({\bar{\mathbb
{P}}})}$, and therefore
$q_m= \sum_{i= 1}^\infty a_ir_i$, where
$a_i$ are the\vspace*{1pt} coefficients $ \langle q_m,r_i \rangle
_{\bar{\mathbb{P}}}$. By
the observation above, we have the equivalent representation $q_m=
\sum_{i= 1}^\infty
\frac{a_i}{\sqrt{\lambda_i}}g_i$.
Therefore the
${L^2(\bar{\mathbb{P}})}$ projection onto
$\mathcal{R}=\operatorname{span}\{r_1,\ldots,r_K\}$ is
${\Pi_{\mathcal{R}}}q_m= \sum_{i= 1}^Ka_i
r_i$, and the ${\bar{\mathcal{H}}}$ projection onto ${\mathcal G}=
\operatorname{span}\{g_1,\ldots,g_K\}$ is ${\Pi_{{\mathcal G} }}q_m=
\sum_{i= 1}^K\frac{a_i}{\sqrt{\lambda_i}}
g_i$. Therefore the relation $g_i=
{\sqrt{\lambda_i}} r_i$ implies
\mbox{$\llVert{\Pi_{\mathcal{R}}} -{\Pi_{{\mathcal G}}} \rrVert_{{\bar
{\mathbb{P}}_n}} = 0$.}
Similar reasoning yields $\llVert{\Pi_{{\tilde{\mathcal
G}}}}q_\ell- {\Pi_{{\tilde{\mathcal V}}}}q_\ell\rrVert
_{{\bar{\mathbb{P}}_n}} = 0$.

\section{Discussion}\label{SecDisc}

In this paper, we have analyzed the performance of spectral clustering
in the context of nonparametric finite mixture models. Our first main
contribution is an upper bound on the distance between the population
level normalized Laplacian embedding and the square root kernelized
density embedding. This bound depends on the maximal similarity
index, the coupling parameter, and the indivisibility parameter.
These parameters all depend on the kernel function, and we present our
analysis for a fixed but arbitrary kernel.

Although this dependence on the kernel function might seem
undesirable, it is actually necessary to guarantee identifiability of
the mixture components in the following sense. A mixture with fully
nonparametric components is a very rich model class: without any
restrictions on the mixture components, any distribution can be
written as a $K$-component mixture in uncountably many ways.
Conversely, when the clustering difficulty function is zero, the
representation of a distribution as a mixture is unique. In
principle, one could optimize over the convex cone of symmetric
positive definite kernel functions so to minimize our clustering
difficulty parameter. In our preliminary numerical experiments, we
have found promising results in using this strategy to choose the
bandwidth in a family of kernels.

Building on our population-level result, we have also provided a result
that characterizes the normalized Laplacian embedding when applied to
a finite collection of $n$ i.i.d. samples. We find that when the
clustering difficulty is small, the embedded samples take on
approximate orthogonal structure: samples from different components
are almost orthogonal with high probability. The emergence of this
form of angular structure allows an angular version of $K$-means
to correctly label most of the samples.

Perhaps surprising is the fact that the optimal bandwidth (minimizing
our upper bound) is nonzero. Although we only provide an upper
bound, we believe this is fundamental to spectral clustering, not an
artifact of our analysis. Again, the principal $K$-dimensional
eigenspace of the Laplacian operator is not a well-defined
mathematical object when the bandwidth is zero. Indeed, as the
bandwidth shrinks to zero, the eigengap distinguishing this eigenspace
from the remaining eigenfucntion vanishes. This eigenspace, however,
is the population-level version of the subspace onto which spectral
clustering projects. For this reason, we caution against shrinking
the bandwidth indefinitely to zero, and we conjecture that there is an
optimal population level bandwidth for spectral clustering. However,
we should mention that we cannot provably rule out the optimality of
an appropriately slowly shrinking bandwidth, and we leave this to
future work. Further investigation of kernel bandwidth selection for
spectral clustering is an interesting avenue for future work.

\section*{Acknowledgments}
The authors thank Johannes Lederer, Elina Robeva, Sivaraman
Balakrishnan, Siqi Wu and Stephen Boyd for helpful discussions. They
are also grateful to the Associate Editor and anonymous referees for
their suggestions that helped improve the manuscript.

\begin{supplement}[id=App]
\sname{Appendix}
\stitle{Remaining proofs and background material\\}
\slink[doi]{10.1214/14-AOS1283SUPP} 
\sdatatype{.pdf}
\sfilename{aos1283\_supp.pdf}
\sdescription{Due to space constraints, we relegate technical details
of the remaining proofs to the supplement~\cite{SchiebingerWainwrightYu14Sup}. This supplementary appendix also gives
an overview of some useful background material, and it includes a
reference list for the symbols used in this paper.}
\end{supplement}



%

\printaddresses
\end{document}